%
\def\conv{\mathop{\vrule height2,6pt depth-2,3pt 
    width 5pt\kern-1pt\rightharpoonup}}

\advance\vsize by 1 true cm
%
\def\dess #1 by #2 (#3){
  \vbox to #2{
    \hrule width #1 height 0pt depth 0pt
    \vfill
    \special{picture #3} 
    }
  }

\def\dessin #1 by #2 (#3 scaled #4){{
  \dimen0=#1 \dimen1=#2
  \divide\dimen0 by 1000 \multiply\dimen0 by #4
  \divide\dimen1 by 1000 \multiply\dimen1 by #4
  \dess \dimen0 by \dimen1 (#3 scaled #4)}
  }
%
\def \trait (#1) (#2) (#3){\vrule width #1pt height #2pt depth #3pt}
\def \fin{\hfill
	\trait (0.1) (5) (0)
	\trait (5) (0.1) (0)
	\kern-5pt
	\trait (5) (5) (-4.9)
	\trait (0.1) (5) (0)
\medskip}
%


\font\sevenbf=cmbx7

\baselineskip=15pt
\abovedisplayskip=15pt plus 4pt minus 9pt
\belowdisplayskip=15pt plus 4pt minus 9pt
\abovedisplayshortskip=3pt plus 4pt
\belowdisplayshortskip=9pt plus 4pt minus 4pt
\let\epsilon=\varepsilon

\def\biblio #1 #2\par{\parindent=30pt\item{}\kern -30pt\rlap{[#1]}\kern
30pt #2\smallskip}
 %
\catcode`\@=11
\def\@lign{\tabskip=0pt\everycr={}}
\def\equations#1{\vcenter{\openup1\jot\displ@y\halign{\hfill\hbox
{$\@lign\displaystyle##$}\hfill\crcr
#1\crcr}}}
\catcode`\@=12
%
\def\pmb#1{\setbox0=\hbox{#1}%
\hbox{\kern-.04em\copy0\kern-\wd0
\kern.08em\copy0\kern-\wd0
\kern-.02em\copy0\kern-\wd0
\kern-.02em\copy0\kern-\wd0
\kern-.02em\box0\kern-\wd0
\kern.02em}}
%
\def\undertilde#1{\setbox0=\hbox{$#1$}
\setbox1=\hbox to \wd0{$\hss\mathchar"0365\hss$}\ht1=0pt\dp1=0pt
\lower\dp0\vbox{\copy0\nointerlineskip\hbox{\lower8pt\copy1}}}
%

%

\def\maj#1#2,{\rm #1\sevenrm #2\rm{}}
\def\Maj#1#2,{\bf #1\sevenbf #2\rm{}}
\outer\def\lemme#1#2 #3. #4\par{\medbreak
\noindent\maj{#1}{#2},\ #3.\enspace{\sl#4}\par
\ifdim\lastskip<\medskipamount\removelastskip\penalty55\medskip\fi}

\def\Remark #1. {\noindent{\Maj REMARK,\ \bf #1. }}

\outer\def\Lemme#1#2 #3. #4\par{\medbreak
\noindent\Maj{#1}{#2},\ \bf #3.\rm\enspace{\sl#4}\par
\ifdim\lastskip<\medskipamount\removelastskip\penalty55\medskip\fi}



\def\Notation #1. {\noindent{\Maj NOTATION,\ \bf #1. }}

\def\Example #1. {\noindent{\Maj EXAMPLE,\ \bf #1. }}

\hfuzz=1cm


\catcode`\ˆ=\active     \def ˆ{\`a}
\catcode`\‰=\active     \def ‰{\^a}
\catcode`\=\active     \def {\c c}
\catcode`\Ž=\active    \def Ž{\'e} 

\catcode`\=\active   \def {\`e}
\catcode`\=\active   \def {\^e}
\catcode`\'=\active   \def '{\"e}
\catcode`\"=\active   \def "{\^\i}
\catcode`\•=\active   \def •{\"\i}
\catcode`\™=\active   \def ™{\^o}
\catcode`\š=\active   \defš{}
\catcode`\=\active   \def {\`u}
\catcode`\ž=\active   \def ž{\^u}
\catcode`\Ÿ=\active   \def Ÿ{\"u}
\catcode`\ =\active   \def  {\tau}
\catcode`\¡=\active   \def ¡{\circ}
\catcode`\¢=\active   \def ¢{\Gamma}
\catcode`\¤=\active   \def ¤{\S\kern 2pt}
\catcode`\¥=\active   \def ¥{\puce}
\catcode`\§=\active   \def §{\beta}
\catcode`\¨=\active   \def ¨{\rho}
\catcode`\©=\active   \def ©{\gamma}
\catcode`\­=\active   \def ­{\neq}
\catcode`\°=\active   \def °{\ifmmode\ldots\else\dots\fi}
\catcode`\±=\active   \def ±{\pm}
\catcode`\²=\active   \def ²{\le}
\catcode`\³=\active   \def ³{\ge}
\catcode`\µ=\active   \def µ{\mu}
\catcode`\¶=\active   \def ¶{\delta}
\catcode`\·=\active   \def ·{\Sigma}
\catcode`\¸=\active   \def ¸{\Pi}
\catcode`\¹=\active   \def ¹{\pi}
\catcode`\»=\active   \def »{\Upsilon}
\catcode`\¾=\active   \def ¾{\alpha}
\catcode`\À=\active   \def À{\cdots}
\catcode`\Â=\active   \def Â{\lambda}
\catcode`\Ã=\active   \def Ã{\sqrt}
\catcode`\Ä=\active   \def Ä{\varphi}
\catcode`\Å=\active   \def Å{\xi}
\catcode`\Æ=\active   \def Æ{\Delta}
\catcode`\Ç=\active   \def Ç{\cup}
\catcode`\È=\active   \def È{\cap}
\catcode`\Ï=\active   \def Ï{\oe}
\catcode`\Ñ=\active   \def Ñ{\to}
\catcode`\Ò=\active   \def Ò{\in}
\catcode`\Ô=\active   \def Ô{\subset}
\catcode`\Õ=\active   \def Õ{\superset}
\catcode`\Ö=\active   \def Ö{\over}
\catcode`\×=\active   \def ×{\nu}
\catcode`\Ù=\active   \def Ù{\Psi}
\catcode`\Ú=\active   \def Ú{\Xi}
\catcode`\Ü=\active   \def Ü{\omega}
\catcode`\Ý=\active   \def Ý{\Omega}
\catcode`\ß=\active   \def ß{\equiv}
\catcode`\à=\active   \def à{\chi}
\catcode`\á=\active   \def á{\Phi}
\catcode`\ä=\active   \def ä{\infty}
\catcode`\å=\active   \def å{\zeta}
\catcode`\æ=\active   \def æ{\varepsilon}
\catcode`\è=\active   \def è{\Lambda}  
\catcode`\é=\active   \def é{\kappa}
\catcode`\ë=\active   \defë{\Theta}
\catcode`\ì=\active   \defì{\eta}
\catcode`\í=\active   \defí{\theta}
\catcode`\î=\active   \defî{\times}
\catcode`\ñ=\active   \defñ{\sigma}
\catcode`\ò=\active   \defò{\psi}

\def\date{\number\day\
\ifcase\month \or janvier \or f\'evrier \or mars \or avril \or mai \or juin \or juillet \or ao\^ut  \or
septembre \or octobre \or novembre \or d\'ecembre \fi
\ \number\year}

\font \Ggras=cmb10 at 12pt
\font \ggras=cmb10 at 11pt

\def\Ga{{\bf \Ggras a}}
\def\Gb{{\bf \Ggras b}}
\def\Gc{{\bf \Ggras c}}
\def\Gd{{\bf \Ggras d}}
\def\Ge{{\bf \Ggras e}}

\def\Gn{{\bf \Ggras n}}

\def\Gt{{\bf \Ggras t}}

\def\Gv{{\bf \Ggras v}}

\def\GA{{\bf \Ggras A}}

\def\GD{{\bf \Ggras D}}
\def\GE{{\bf \Ggras E}}

\def\GI{{\bf \Ggras I}}

\def\GL{{\bf \Ggras L}}

\def\GS{{\bf \Ggras S}}

\def\GU{{\bf \Ggras U}}
\def\GV{{\bf \Ggras V}}

\def\GR{{\bf \Ggras R}}

\def\liminf{\mathop{\underline{\rm lim}}}
\def\limsup{\mathop{\overline{\hbox{\rm lim}}}}
\def\sym{\fam\comfam\com}
\font\tensym=msbm10
\font\sevensym=msbm7
\font\fivesym=msbm5
\newfam\symfam
\textfont\symfam=\tensym
\scriptfont\symfam=\sevensym
\scriptfont\symfam=\fivesym
\def\sym{\fam\symfam\relax}
\def\N{{\sym N}}
\def\R{{\sym R}}

\def\D{{\sym D}}
\def\E{{\sym E}}
\let\ds\displaystyle

\def\liminf{\mathop{\underline{\rm lim}}}
\def\limsup{\mathop{\overline{\hbox{\rm lim}}}}
\def\sym{\fam\comfam\com}
\font\tensym=msbm10
\font\sevensym=msbm7
\font\fivesym=msbm5
\newfam\symfam
\textfont\symfam=\tensym
\scriptfont\symfam=\sevensym
\scriptfont\symfam=\fivesym
\def\sym{\fam\symfam\relax}
\def\N{{\sym N}}
\def\R{{\sym R}}

\def\D{{\sym D}}
\def\E{{\sym E}}
\let\ds\displaystyle

\centerline{\Ggras A simplified model for elastic  thin shells.}
\bigskip
\bigskip
\centerline{\it Dominique Blanchard ${ }^a$ and Georges Griso ${ }^b$}

\centerline{\it ${ }^a$Universit\'e de Rouen, UMR 6085,  76801   Saint Etienne du Rouvray Cedex, France,}

\centerline{\it e-mail: dominique.blanchard@univ-rouen.fr,  blanchar@ann.jussieu.fr}

\centerline{\it ${ }^b$ Laboratoire J.-L. Lions--CNRS, Bo\^\i te courrier 187, Universit\'e  Pierre et
Marie Curie,}

\centerline{\it  4 place Jussieu, 75005 Paris, France, e-mail: griso@ann.jussieu.fr}

\bigskip
\noindent{\Ggras Abstract. } We introduce a simplified model for the minimization of the elastic energy in thin shells. This model is not obtained by an  asymptotic analysis.  The nonlinear simplified model admits always minimizers by contrast with the original one. We show the relevance of our approach by proving that the rescaled minimum of the simplified model and the rescaled infimum of the full model have the same limit as the thickness tends to $0$. The simplified energy can be expressed as a functional acting over   fields defined on the mid-surface of the shell and where  the thickness remains as a parameter.

\bigskip
\noindent {\ggras Keywords: } {nonlinear elasticity, shells.}

\noindent{\ggras  2000 MSC: } 74B20, 74K20, 74K25, 49J45.

\bigskip

\noindent {\Ggras 1.  Introduction}\medskip
\medskip

This paper is devoted to introduce and  justify a simplified model for nonlinear elastic shells.  Let $\omega$ be a bounded Lipschitz domain of $\R^2$ and $\phi$ be a smooth function from $\overline{\omega}$ into $\R^3$ (see the detailed assumptions on $\phi$ in Section 2) and set $S=\phi(\omega)$.  We denote by $\Gn$ an unit vector field normal to $S$ and by $\Phi$ the map $(s_1,s_2,s_3)\longrightarrow \phi(s_1,s_2)+s_3\Gn(s_1,s_2)$. The elastic shell is defined by ${\cal Q}_\delta=\Phi(\omega \times ]-\delta,\delta[)$ and we consider  that it is  clamped on a part of its lateral boundary $\Gamma_{0,\delta}=\Phi(\gamma_0\times ]-\delta,\delta[)$, where $\gamma_0\subset \partial \omega$. The  energy density is denoted $W$ and we assume that ${\cal Q}_\delta$ is  submitted to applied body forces $f_{\kappa,\delta}$ whose order with respect to $\delta$ depends upon a parameter $\kappa$ (see the order of $f_{\kappa,\delta}$ below). The total energy  is given by $\ds J_{\kappa,\delta}(v)=\int_{{\cal Q}_\delta} W(E( v))-\int_{{\cal Q}_\delta} f_{\kappa,\delta}\cdot (v-I_d)$ if $\det (\nabla v)>0$ and where $E(v)={1/2}\big((\nabla v)^T\nabla v-\GI_3\big)$ is   the Green-St Venant's tensor and $I_d$ is the identity map. We set
$$ m_{\kappa,\delta}=\inf_{v\in \GU_\delta}J_{\kappa,\delta}(v),
$$
where $\GU_\delta$ is the set of admissible deformations (which are equal to the identity map on $\Gamma_{0,\delta}$). The Korn's type inequalities  established in [6] (see also [12]) allow us to prove that if the order of  $f_{\kappa,\delta}$ is equal to $\delta^{2\kappa-2}$ for $1\le\kappa\le2$ (or $\delta^{\kappa}$ for $\kappa\ge 2$), then  the  order of $m_{\kappa,\delta}$ is   $\delta^{2\kappa-1}$. 

 Even for a classical St-Venant-Kirchhoff's material, proving the existence of a minimizer for $J_{\kappa,\delta}$  is still an open problem. The aim of this paper is to replace the above minimization problem  by a minimization problem for a simplified functional $J^s_{\kappa,\delta}$  defined on a  new set  $\D_{\delta,\gamma_0}$  and which admits  a minimum
 $$ m^s_{\kappa,\delta}=\min_{\Gv\in \D_{\delta,\gamma_0}} J^s_{\kappa,\delta}(\Gv)$$ of the same  order as $m_{\kappa,\delta}$. This approximation is justified if one shows  that 
 $$\lim_{\delta\to 0}{m_{\kappa,\delta}-m^s_{\kappa,\delta} \over \delta^{2\kappa-1}}=0.$$
 In the present paper we show this result in the case $\kappa=2$ (other critical cases will be investigated in forthcoming papers).
 
 The expression of $J^s_{\kappa,\delta}$ and the choice of $\D_{\delta,\gamma_0}$ rely on the decomposition technique introduced in [6].  Let us  recall that  a deformation $v$ of the  shell ${\cal Q}_\delta$, whose "geometrical energy" -$||dist(\nabla v,SO(3))||_{L^2({\cal Q}_\delta)}$- is at most of order $\delta^{3/2}$,  is decomposed as (see [6] or Theorem 3.1 below)
$$ v(x)={\cal V}(s_1,s_2) + s_3\GR(s_1,s_2) \Gn(s_1,s_2)+\overline{v}(s_1,s_2,s_3),\quad x=\Phi(s),\quad \hbox{for a.e. } s=(s_1,s_2,s_3)\in \omega\times]-\delta,\delta[.$$ The field ${\cal V}$ stands for the mid-surface deformation, the  matrix field $\GR$ takes its values in $SO(3)$ and represents the rotations of the fibers and $\overline{v}$ is the warping of theses fibers. It is also shown in [6] that the fields ${\cal V}$, $\GR$ and  $\overline{v}$ satisfy the natural boundary conditions on $\gamma_0$ and on $\gamma_0\times ]-\delta,\delta[$ and that they are  estimated in terms of  $||dist(\nabla v,SO(3))||_{L^2({\cal Q}_\delta)}$ and $\delta$.  With the help of these estimates, we justify the simplification of the Green-St Venant's strain tensor $E(v)$ in order to give a simplified matrix $\widehat{E}(\Gv)$ which depends on the triplet $\Gv=({\cal V},\GR,\overline{v})$ associated to a deformation $v$. This matrix depends linearly upon $\ds{\partial\overline{v}\over \partial s_3}$ and on the first partial derivatives of ${\cal V}$ and $\GR$ (see Section 5) but  which is nonlinear with respect to $({\cal V},\GR,\overline{v})$. 

Then we define the set $\D_{\delta,\gamma_0}$  of  admissible  triplets $\Gv=({\cal V},\GR,\overline{v})$ and we derive the simplified total energy $J^s_{\kappa,\delta}(\Gv)$ as follows. 
Firstly we replace $\int_{{\cal Q}_\delta} W(E( v))$ by $\int_{{\cal Q}_\delta} Q(E( \Gv))$ where $Q$ is a quadratic form which is assumed to approximate $W$ near the origin. Secondly we add two penalization terms in order to approach  the usual limit kinematic condition $\ds{\partial {\cal V}\over \partial s_\alpha}=\GR\Gt_\alpha$ and to insure the coerciveness of  $J^s_{\kappa,\delta}$. Finally in the term involving the forces we neglect the contribution of the warping $\overline{v}$. As announced above we prove that  $J^s_{\kappa,\delta}$ admit minimizers on $\D_{\delta,\gamma_0}$.  We justify the approximation process described above in the case  $\kappa=2$. 
\smallskip
In some sense, the introduction of $J^s_{\kappa,\delta}$ can be seen as  a nonlinear version of the approach which leads to  the simplified Timoshenko's model for rods, the  Reisner-Mindlin's model  for plates and  the Koiter's model for shells in linear elasticity.

As general references on the theory of nonlinear elasticity, we refer to  [8] and [24] and to the extensive bibliographies of these works. A general theory for the existence of minimizers of nonlinear elastic energies can be found in [1]. For the justification  of  plate or shell models in nonlinear elasticity we refer  to [9], [10], [11], [13], [15], [18], [23], [25], [26]. The derivation of limit energies for thin  domains using $\Gamma$-convergence arguments are developed in [14], [15], [22], [23]. The decomposition of the deformations in thin structures is introduced in [17], [18] and  a few applications to the junctions of multi-structures and homogenization are given in  [2], [3], [4]. The justification of simplified models for rods and plates in linear elasticity, based on a decomposition technique of the displacement, is presented in [19], [20]. In this linear case, error estimates between the solution of the initial model and the one of the simplified model  are also established. In some sense, these works give a mathematical justification of  Timoshenko's model for rods and Reisner-Mindlin's model  for plates.
 
 The paper is organized as follows. Section 2 is devoted to describe the geometry of the shell and to give a few notations.  In Section 3 we recall the results of [6]:   decomposition of a deformation of a thin shell,  estimates on the terms of this decomposition and  two nonlinear Korn's type inequalities. Section 4 is concerned with a standard rescaling. We present the simplification  of the Green-St Venant's strain tensor of a deformation in Section 5. We  also introduce  the set $\D_{\delta,\gamma_0}$ of admissible triplets $\Gv=({\cal V},\GR,\overline{V})$ and we  prove Korn's type inequalities for the elements of $\D_{\delta,\gamma_0}$ (see Corollary 5.3). In Section 6 we consider nonlinear elastic shells and we use the results of [6]  to scale the applied forces in order to obtain a priori estimates on $m_{\kappa,\delta}$.  Section 7 is devoted to introduce the simplified  energy $J^s_{\kappa,\delta}$ and to prove the existence of minimizers. In Sections 8 and 9, we restrict the analysis to $\kappa=2$. We prove that  
 $$\lim_{\delta\to 0}{m_{2,\delta} \over \delta^{3}}=\lim_{\delta\to 0}{m^s_{2,\delta} \over \delta^{3}}=m^s_2$$ where $m^s_2$ is the minimum of a functional defined over a set of triplets. 
 In Section 10, we give an alternative formulation of the minimization problem for $J^s_{\kappa,\delta}$ through elimination of the variable ${\overline v}$. Then we obtain that $m^s_{\kappa,\delta}$ is the minimum of a functional which depends only upon $({\cal V}, \GR)$. At last  an appendix contains  an approximation result for the elements of $\D_{\delta,\gamma_0}$ and  an algebraic elimination process for quadratic forms.
The results of this paper were announced in [7].
\bigskip
\noindent {\Ggras 2. The geometry and notations.}
\medskip
Let us introduce a few  notations and definitions concerning the geometry of the shell.

 Let $\omega$ be a bounded domain in $\R^2$ with lipschitzian boundary and let $\phi$ be an injective mapping   from $\overline{\omega}$ into $\R^3$ of class ${\cal C}^2$.  We denote $S$ the surface $\phi(\overline{\omega})$.  We assume that the two vectors $\displaystyle {\partial \phi\over \partial s_1}(s_1,s_2)$ and $\displaystyle {\partial \phi\over \partial s_2}(s_1,s_2)$ are linearly independent at each point  $(s_1,s_2)\in\overline{\omega}$. 
\smallskip
We set
$$\Gt_1={\partial\phi\over \partial s_1},\qquad \Gt_2={\partial\phi\over \partial s_2},\qquad \Gn={\displaystyle
\Gt_1\land\Gt_2\over \displaystyle\bigl\| \Gt_1\land \Gt_2\bigr\|_2}.\leqno(2.1)$$  The vectors $\Gt_1$ and  $\Gt_2$ are tangential vectors to the surface $S$ and the vector $\Gn$ is a unit normal vector  to this surface.  We set
 $$ \Omega_\delta=\omega\times ]-\delta,\delta[.$$

 Now we consider the mapping $\Phi\, : \,\overline{\omega}\times \R  \longrightarrow \R^3$ defined by
$$\Phi\enskip:\enskip (s_1,s_2,s_3) \longmapsto x=\phi(s_1,s_2)+s_3\Gn(s_1,s_2).\leqno(2.2)$$

There exists  $\delta_0\in (0,1]$  depending only on $S$, such that the restriction of  $\Phi$ to the compact set
$\overline{\Omega}_{\delta_0}=\overline{\omega}\times [-\delta_0,\delta_0]$ is a ${\cal C}^1-$ diffeomorphism of that set onto its range (see e.g. [21]). Hence, there exist two constants $c_0>0$ and $c_1\ge c_0$, which depend only on $\phi$, such that
$$ \forall s\in  \Omega_{\delta_0},\quad c_0\le |||\nabla_s\Phi(s))||| \le c_1,\quad \hbox{and for } x=\Phi(s) \qquad c_0\le |||\nabla_x\Phi^{-1}(x))||| \le c_1.\leqno(2.3)$$  
\smallskip 
\noindent{\Ggras Definition 2.1. }{\it For $ \delta\in (0, 
\delta_0]$, the shell  ${\cal Q}_\delta$ is defined as follows:
$${\cal Q}_\delta=\Phi(\Omega_\delta).$$} The mid-surface of the shell is $S$. The  fibers of the shell are the segments $\Phi\big(\{(s_1,s_2)\}\times]-\delta,\delta[\big)$, $(s_1,s_2)\in \omega$. The lateral boundary of  the  shell is $\Gamma_\delta=\Phi(\partial\omega\times]-\delta,\delta[)$. In the following sections the shell will be fixed on a part of its lateral boundary. Let $\gamma_0$ be an open subset of $\partial \omega$ which made of a finite number of connected components (whose closure are disjoint). We assume that the shell is clamped on 
 $$\Gamma_{0,\delta}=\Phi(\gamma_0\times ]-\delta,\delta[).$$ The admissible deformations $v$ of the shell must then satisfy
$$v=I_d\qquad \hbox{on}\quad \Gamma_{0,\delta}\leqno(2.4)$$ where  $I_d$ is  the identity map of $\R^3$.
\medskip
\noindent{\Ggras Notation. } From now on we denote by $c$ and $C$ two positive generic constants which do not depend on $\delta$.
We respectively  note by $x$ and $s$ the generic points of ${\cal Q}_\delta$ and of $\Omega_\delta$.  A field $v$ defined on ${\cal Q}_\delta$ can be also considered as  a field defined on $\Omega_\delta$ that, as a convention,   we will also denote by $v$. As far as the gradients of field $v$, say in $(W^{1,1}({\cal Q}_\delta))^3$,  are concerned we have $\nabla_x v$ and $\nabla_s v=\nabla_x v .\nabla \Phi$  for a.e. $x=\Phi(s)$ and (2.3) shows that 
$$ c|||\nabla_x v(x) |||\le |||\nabla_s v(s)|||\le C|||\nabla_x v(x) |||.$$
\medskip
\noindent {\Ggras 3. Korn's type inequalities for shells. Decomposition of a deformation.}
\smallskip
We first recall the  Korn's type  inequalities for shells established in Section 4 of [6]. Let $v$ be an  admissible deformation  belonging to $\big(H^1({\cal Q}_\delta) \big)^3$ and satisfying the boundary condition (2.4). Setting  $\ds {\cal V}(s_1,s_2)={1\over 2\delta} \int_{-\delta}^{\delta} v(s_1,s_2, t)dt\;$  a.e. $(s_1,s_2)\in\omega$, we have 
$$\left\{\eqalign{ 
 &||v-I_d||_{(L^2({\cal Q}_\delta))^3}+||\nabla_x v-\GI_3||_{(L^2({\cal Q}_\delta))^9}\le C\big( \delta^{1/2}+||dist(\nabla_x v,SO(3))||_{L^2({\cal Q}_{\delta})}\big),\cr
 & ||(v-I_d)-({\cal V}-\phi)||_{(L^2({\cal Q}_\delta))^3}\le  C\delta\big( \delta^{1/2}+||dist(\nabla_x v,SO(3))||_{L^2({\cal Q}_{\delta})}\big), \cr}\right.\leqno (3.1)$$ and
 $$\left\{\eqalign{ 
 & ||v-I_d||_{(L^2({\cal Q}_\delta))^3}+||\nabla_x v-\GI_3||_{(L^2({\cal Q}_\delta))^9}\le {C\over \delta}||dist(\nabla_x v,SO(3))||_{L^2({\cal Q}_{\delta})},\cr
& ||(v-I_d)-({\cal V}-\phi)||_{(L^2({\cal Q}_\delta))^3}\le  C ||\hbox{dist}(\nabla_x v,SO(3))||_{L^2({\cal Q}_\delta)}.\cr}\right.\leqno (3.2)$$ Inequalities (3.1) are better than  those (3.2) if the order of the geometric energy $||dist(\nabla_x v,SO(3))||_{L^2({\cal Q}_{\delta})}$ is greater than $\delta^{3/2}$.

 Now the  theorem of decomposition of the deformations established in [6] (see Theorem 3.4 of Section 3) is given below.
\medskip
\noindent{\ggras  Theorem 3.1.  }{\it There exists a constant $C(S)$ which depends only on the mid-surface of the shell such that for all deformation   $v$ belonging to $\big(H^1({\cal Q}_\delta) \big)^3$ and satisfying$$||dist(\nabla_xv ,SO(3))||_{L^2({\cal Q}_\delta)}\le C(S)\delta^{3/2}, \leqno(3.3)$$ 
then, there exist ${\cal V}\in (H^1(\omega))^3$, $\GR\in\big(H^1(\omega) \big)^{3\times 3}$ satisfying $\GR(s_1,s_2)\in SO(3)$  for a.e. $(s_1,s_2)\in \omega$ and $\overline{v}$ belonging to $ \big(H^1({\cal Q}_\delta) \big)^3$   such that for a.e. $s\in \Omega_\delta$ 
 $$ v(s)={\cal V} (s_1,s_2)+ s_3\GR(s_1,s_2) \Gn(s_1,s_2)+\overline{v}(s),\leqno (3.4)$$
where we  can choose $\ds {\cal V}(s_1,s_2)={1\over 2\delta} \int_{-\delta}^{\delta} v(s_1,s_2, t)dt\;$  a.e. $(s_1,s_2)\in\omega$,
and such that the following estimates hold:
$$\left\{\eqalign{
&||\overline{v}||_{(L^2(\Omega_\delta))^3}\le C\delta ||dist(\nabla_x v,SO(3))||_{L^2({\cal Q}_\delta)}\cr
&||\nabla_s\overline{v}||_{(L^2(\Omega_\delta))^9}\le C ||dist(\nabla_x v,SO(3))||_{L^2({\cal Q}_\delta)}\cr
&\Bigl\|{\partial \GR\over \partial s_\alpha}\Big\|_{(L^2(\omega))^9}\le {C\over \delta^{3/2}} ||dist(\nabla_x v,SO(3))||_{L^2({\cal Q}_\delta)}\cr
& \Bigl\|{\partial{\cal V}\over \partial s_\alpha}-\GR \Gt_\alpha\Big\|_{(L^2(\omega))^3}\le {C\over \delta^{1/2}}||dist(\nabla_x v,SO(3))||_{L^2({\cal Q}_\delta)}\cr
& \bigl\|\nabla_x v-\GR \big\|_{(L^2(\Omega_\delta))^9}\le C||dist(\nabla_x v,SO(3))||_{L^2({\cal Q}_\delta)}.\cr}\right.\leqno(3.5)$$\fin }
Due to (3.4) and to the definition of ${\cal V}$, the field $\overline {v}$ satisfies $\ds  \int_{-\delta}^{\delta} \overline{v}(s_1,s_2, t)dt=0\;$  a.e. $(s_1,s_2)\in\omega$.

\medskip
If the deformation $v$ as in Theorem 3.1 satisfies the boundary condition (2.4) then indeed
 $$ {\cal V}=\phi\qquad \hbox{on}\quad \gamma_0.\leqno(3.6)$$

 Moreover due to Lemma 4.1 of [6], we can choose $\GR$ and $\overline{v}$ in Theorem 3.1 above such that
 $$\GR=\GI_3 \qquad \hbox{on } \gamma_0, \qquad \overline{v}=0 \qquad \hbox{on }  \Gamma_{0,\delta}.\leqno (3.7)$$

\noindent From estimates (3.5) we also derive the following ones:
$$\left\{\eqalign{
||\GR-\GI_3||_{(L^2(\omega))^9}&\le {C\over \delta^{3/2}} ||dist(\nabla_x v,SO(3))||_{L^2({\cal Q}_\delta)}\cr
||{\cal V}-\phi||_{(L^2(\omega))^3}&\le {C\over \delta^{3/2}} ||dist(\nabla_x v,SO(3))||_{L^2({\cal Q}_\delta)}.\cr}\right.\leqno(3.8)$$ 



%

\noindent {\Ggras 4.  Rescaling $\Omega_\delta$ }
  \medskip
As usual, we rescale $\Omega_\delta$ using the operator 
$$(\Pi_\delta w)(s_1,s_2,S_3)=w(s_1,s_2,\delta S_3)\hbox{ for any}\;\; (s_1,s_2,S_3)\in \Omega $$
defined for e.g. $w\in L^2(\Omega_\delta)$ for which $ (\Pi_\delta w)\in  L^2(\Omega)$. Let $v$ be a deformation decomposed as (3.4), by transforming by $\Pi_\delta$  we obtain
$$\Pi_\delta (v)(s_1,s_2,S_3)={\cal V} (s_1,s_2)+ \delta S_3\GR(s_1,s_2) \Gn(s_1,s_2)+\Pi_\delta (\overline{v})(s_1,s_2,S_3),\qquad \hbox{for a.e. }(s_1,s_2,S_3)\in \Omega.$$
The estimates (3.5) of  $\overline{v}$  transposed over $\Omega$ are (notice that $\ds \Pi_\delta\Big({\partial \overline{v}\over \partial s_3}\Big)={1\over \delta}{\partial\Pi_\delta( \overline{v})\over \partial S_3}$)
$$\left\{\eqalign{
&||\Pi_\delta (\overline{v})||_{(L^2(\Omega))^3}\le C \delta^{1/2}||dist(\nabla_x v,SO(3))||_{L^2({\cal Q}_\delta)}\cr
&\Big\|{\partial \Pi_\delta (\overline{v})\over \partial s_1}\Big\|_{(L^2(\Omega))^3}\le {C\over \delta^{1/2}} ||dist(\nabla_x v,SO(3))||_{L^2({\cal Q}_\delta)}\cr 
&\Big\|{\partial \Pi_\delta (\overline{v})\over \partial s_2}\Big\|_{(L^2(\Omega))^3}\le  {C\over \delta^{1/2}}  ||dist(\nabla_x v,SO(3))||_{L^2({\cal Q}_\delta)} \cr 
&\Big\|{\partial \Pi_\delta (\overline{v})\over \partial S_3}\Big\|_{(L^2(\Omega))^3}\le C \delta^{1/2} ||dist(\nabla_x v,SO(3))||_{L^2({\cal Q}_\delta)}, \cr}\right.\leqno(4.1)$$
\noindent {\Ggras  5.  Simplification in the Green-St Venant's strain tensor. }
  \medskip
\noindent In this section we introduce a simplification of the  Green-St Venant's strain tensor $E(v)={1/2}\big((\nabla v)^T\nabla v-\GI_3\big)$. Let $v$ be a deformation of the shell  belonging to $(H^1({\cal Q}_\delta))^3$ and satisfying the condition (3.3). We decompose $v$ as  (3.4). We have the identity
$$(\nabla_x v)^T\nabla_x v-\GI_3=(\nabla_x v-\GR)^T\GR+\GR^T(\nabla_x v-\GR)+(\nabla_x v-\GR)^T(\nabla_x v-\GR).$$
In order to compare the orders (of the norms)  of the different terms in the above equality, we work in the fix domain $\Omega$ using the operator $\Pi_\delta$. Thanks to estimates (3.5) we get
$$\left\{\eqalign{
||\Pi_\delta\big((\nabla_x v-\GR)^T\GR+\GR^T(\nabla_x v-\GR)\big)||_{(L^1(\Omega))^{3\times 3}}\le C{||dist(\nabla_x v,SO(3))||_{L^2({\cal Q}_\delta)}\over \delta^{1/2}},\cr
||\Pi_\delta\big((\nabla_x v-\GR)^T(\nabla_x v-\GR)\big)||_{(L^1(\Omega))^{3\times 3}}\le C\Big[{||dist(\nabla_x v,SO(3))||_{L^2({\cal Q}_\delta)}\over \delta^{1/2}}\Big]^2.\cr}\right.\leqno(5.1)$$ In view of (3.1), these estimates show that the term $\Pi_\delta\big((\nabla_x v-\GR)^T(\nabla_x v-\GR)\big)$ can be neglected in $E(v)$.  

Now we have
$${\partial v\over \partial s_1}=\nabla_x v \,\big(\Gt_1+s_3{\partial \Gn\over \partial s_1}\big), \qquad  {\partial v\over\partial s_2}=\nabla_x v\,\big(\Gt_2+s_3{\partial \Gn\over \partial s_2}\big),
\qquad  {\partial v\over\partial s_3}=\nabla_x v\,\Gn.$$  Then
$$\Pi_\delta(\nabla_x v-\GR) \,\big(\Gt_\alpha +\delta S_3{\partial \Gn\over \partial s_\alpha }\big)=\Big({\partial{\cal V}\over \partial s_\alpha }-\GR\Gt_\alpha\Big)+\delta S_3{\partial \GR\over \partial s_\alpha}\Gn+  {\partial\Pi_\delta\overline{v}\over \partial s_\alpha}, \qquad  \Pi_\delta (\nabla_x v-\GR)\,\Gn= {1\over \delta}{\partial \Pi_\delta \overline{v}\over\partial S_3}.$$ First, we can neglect the term $\ds\delta S_3{\partial \Gn\over \partial s_\alpha}$ which is of order $\delta$ in the quantity $\ds\Gt_\alpha +\delta S_3{\partial \Gn\over \partial s_\alpha }$. Secondly, as a  consequence of these equalities and the following estimates (obtained from (3.5) and (4.1)):
$$\left\{\eqalign{
\Big\|{\partial{\cal V}\over \partial s_\alpha }-\GR\Gt_\alpha\Big\|_{(L^2(\Omega))^3}\le C{||dist(\nabla_x v,SO(3))||_{L^2({\cal Q}_\delta)}\over \delta^{1/2}},\cr
\Big\| \delta S_3{\partial \GR\over \partial s_\alpha}\Gn\Big\|_{(L^2(\Omega))^3}\le C{||dist(\nabla_x v,SO(3))||_{L^2({\cal Q}_\delta)}\over \delta^{1/2}},\cr
\Big\|{1\over \delta}{\partial \Pi_\delta\overline{v}\over \partial S_3}\Big\|_{(L^2(\Omega))^3}\le C{||dist(\nabla_x v,SO(3))||_{L^2({\cal Q}_\delta)}\over \delta^{1/2}},\cr
\Big\|{\partial \Pi_\delta\overline{v}\over \partial s_\alpha} \|_{(H^{-1}(\Omega))^3}\le C \delta\Big[{||dist(\nabla_x v,SO(3))||_{L^2({\cal Q}_\delta)}\over \delta^{1/2}}\Big]\cr}\right.\leqno(5.2)$$ we deduce that in the quantity $\Pi_\delta\big((\nabla_x v-\GR)^T\GR+\GR^T(\nabla_x v-\GR)\big)$ we can neglect  the terms  $\ds {\partial \Pi_\delta \overline{v}\over \partial s_\alpha}$.
 \smallskip

 Now, if in  the Green-St Venant's strain tensor of  $v$  we carry out the simplifications mentioned above, we are brought to replace  
$$\eqalign{
&{1\over 2 }\Pi_\delta\big((\nabla_x v )^T\nabla_x v-\GI_3\big)\quad \hbox{ by }\quad (\Gt_1\,|\, \Gt_2\,|\, \Gn)^{-T}\Pi_\delta\big(E^s(v)\big)(\Gt_1\,|\, \Gt_2\,|\, \Gn)^{-1}\cr
\hbox{or } \quad &{1\over 2 } \big((\nabla_x v )^T\nabla_x v-\GI_3\big)\quad \hbox{ by }\quad (\Gt_1\,|\, \Gt_2\,|\, \Gn)^{-T} E^s(v)(\Gt_1\,|\, \Gt_2\,|\, \Gn)^{-1}\cr}$$ where the symmetric matrix $E^s(v)\in (L^2(\Omega_\delta))^{3\times 3}$ is equal to
$$\eqalign{
E^s(v)&=\pmatrix{
\displaystyle  s_3\Gamma_{11}(\GR)+{\cal Z}_{11} & \displaystyle   s_3\Gamma_{12}(\GR)+ {\cal Z}_{12} &\displaystyle  {1\over 2}\GR^T{\partial\overline{v} \over \partial s_3} \cdot \Gt_1+{1\over 2 }{\cal Z}_{31} \cr
 * & \displaystyle  s_3\Gamma_{22}(\GR)+ {\cal Z}_{22}   &\displaystyle  {1\over 2}\GR^T{\partial\overline{v} \over \partial s_3} \cdot\Gt_2+{1\over 2 }{\cal Z}_{32}\cr
 * & *&  \displaystyle \GR^T{\partial\overline{v} \over \partial  s_3} \cdot\Gn \cr}\cr
\Gamma_{\alpha\beta}(\GR)&={1\over 2} \Big[{\partial \GR\over \partial s_\alpha} \Gn \cdot\GR\Gt_\beta+{\partial \GR\over \partial s_\beta} \Gn \cdot\GR\Gt_\alpha\Big],\cr
{\cal Z}_{\alpha\beta}&={1\over 2} \Big[\Big({\partial{\cal V}\over \partial s_\alpha }-\GR\Gt_\alpha\Big)\cdot\GR\Gt_\beta+\Big({\partial{\cal V}\over \partial s_\beta}-\GR\Gt_\beta\Big)\cdot\GR\Gt_\alpha\Big],\quad  {\cal Z}_{3\alpha}= {\partial{\cal V}\over \partial s_\alpha }\cdot\GR\Gn,\cr}\leqno(5.3)$$
  where $(\Gt_1\,|\, \Gt_2\,|\, \Gn)$ denotes the $3\times 3$ matrix with first column   $\Gt_1$, second column $\Gt_2$ and  third column $\Gn$ and  where $(\Gt_1\,|\, \Gt_2\,|\, \Gn)^{-T}=\big((\Gt_1\,|\, \Gt_2\,|\, \Gn)^{-1}\big)^{T}$. Let us notice that $E^s(v)$ belongs to $\big(L^2(\Omega_\delta)\big)^{3\times 3}$ for any deformation 
$$ v(s)={\cal V}(s_1,s_2)+s_3\GR(s_1,s_2)\Gn(s_1,s_2)+\overline{v}(s),\qquad \hbox{for a.e. } s\in \Omega_\delta$$  where ${\cal V}\in (H^1(\omega))^3$, $\GR\in H^1(\omega ; SO(3))$ and $\overline{v}\in L^2(\omega ; H^1(-\delta , \delta))^3$.
\medskip
\noindent{\ggras Remark 5.1.}  From the   last estimate in (3.5) we deduce that
$$||\Pi_\delta\big(\nabla_x v-\GR\big)||_{(L^2(\Omega))^{3\times 3}}\le C\delta $$ and then we get that the set
$$\bigl\{s\in \Omega\; |\; |||\Pi_\delta\big(\nabla_x v-\GR\big)(s)||| \ge 1\big\}$$ has  a measure less than $C\delta^2$.  It follows that   the measure of the set 
$$\bigl\{s\in \Omega\; |\;  \det\big(\Pi_\delta\big[\nabla_x v\big](s)\big)\le 0 \big\}$$ tends to $0$ as $\delta$ goes to $0$. \fin
\medskip
Now, we introduce the following closed  subset $\D_\delta$ of $ (H^1(\omega))^3\times (H^1(\omega))^{3\times 3}\times ( L^2(\omega ; H^1(-\delta,\delta)))^3$
$$\eqalign{
\D_\delta=\Big\{ & \Gv=({\cal V} , \GR , \overline{v})\in  (H^1(\omega))^3\times (H^1(\omega))^{3\times 3}\times ( L^2(\omega ; H^1(-\delta,\delta)))^3\; |\;\cr
& \quad \GR(s_1, s_2)\in SO(3),\qquad  \int_{-\delta}^\delta\overline{v}(s_1,s_2,s_3)ds_3=0,\cr
&  \int_{-\delta}^\delta s_3\overline{v}(s_1,s_2,s_3)\cdot\Gt_\alpha(s_1,s_2) ds_3=0,\enskip  \hbox{for a.e. } (s_1,s_2)\in \omega, \enskip \alpha=1,2.\Big\}\cr}$$  
The last condition on $\overline{v}$ in $\D_\delta$ is not satisfied in general (if $\overline{v}$ is the warping introduced in Theorem 3.1), loosely speaking  this new condition will allow to decouple the estimates of $\overline{v}$ and ${\cal Z}_{i\alpha}$ (see the proof of  Proposition 5.2).

For any $\Gv\in \D_\delta$, we consider  $v$  defined by 
$$\eqalign{
v(s)={\cal V}(s_1,s_2)+s_3\GR(s_1,s_2)\Gn(s_1,s_2)+\overline{v}(s),\cr}\qquad \hbox{for a.e. } s\in \Omega_\delta.\leqno(5.4)$$ 
The deformation $v$ belongs to $(L^2(\omega ; H^1(-\delta, \delta)))^3$ so that, in general, the Green-St Venant's tensor of $v$ is not defined. Nevertheless, the tensor field $E^s(v)$ belongs to $(L^2(\Omega_\delta))^{3\times 3}$ and we set 
$$\widehat{E}(\Gv)=E^s(v),\qquad \qquad \widehat{E}(\Gv)\in (L^2(\Omega_\delta)^{3\times 3}.\leqno(5.5)$$ 
\smallskip
Let us point out that if a triplet $\Gv$ satisfies the limit kinematic condition $\ds {\partial{\cal V}\over \partial s_\alpha }=\GR\Gt_\alpha$, then it is easy to obtain
$${1\over \delta}||\overline{v}||_{(L^2(\Omega_\delta)^3}+ \Big\|{\partial\overline{v}\over \partial s_3}\Big\|_{(L^2(\Omega_\delta)^3}\le  ||\widehat{E}(\Gv)||_{(L^2(\Omega_\delta))^{3\times 3}},\qquad \Big\|{\partial \GR\over \partial s_\alpha}\Big\|_{(L^2(\omega))^{3\times 3}}
 \le {C\over \delta^{3/2}}||\widehat{E}(\Gv)||_{(L^2(\Omega_\delta))^{3\times 3}}$$ which permits with some boundary conditions to control  the product norm of $\Gv$ in term of $||\widehat{E}(\Gv)||_{(L^2(\Omega_\delta))^{3\times 3}}$ and $\delta$. In order to define an energy which have this property for any $\Gv\in \D_\delta$,  we are led to add two penalization terms, which vanish as $\delta \to 0$, to $|| \widehat{E}(\Gv)||^2_{(L^2(\Omega_\delta))^{3\times 3}}$.  This is why for every deformation $\Gv\in \D_\delta$ we set
$${\cal E}_\delta(\Gv)=||\widehat{E}(\Gv)||^2_{(L^2(\Omega_\delta))^{3\times 3}}+  \delta^3\Big\|{\partial \GR\over \partial s_1} \Gt_2-{\partial \GR\over \partial s_2} \Gt_1\Big\|^2_{(L^2(\omega))^{3}}+\delta \Big\|{\partial{\cal V}\over \partial s_1 }\cdot\GR\Gt_2- {\partial{\cal V}\over \partial s_2}\cdot\GR\Gt_1\Big\|^2_{L^2(\omega)}.\leqno (5.6)$$ 
\noindent {\Ggras Proposition 5.2. } {\it There exists a positive constant $C$ which does not depend on $\delta$ such that for all $\Gv\in \D_\delta$
$$\eqalign{
{1\over \delta}||\overline{v}||_{(L^2(\Omega_\delta)^3}+ \Big\|{\partial\overline{v}\over \partial s_3}\Big\|_{(L^2(\Omega_\delta)^3}\le  ||\widehat{E}(\Gv)||_{(L^2(\Omega_\delta))^{3\times 3}}\cr
\Big\|{\partial \GR\over \partial s_\alpha}\Big\|^2_{(L^2(\omega))^{3\times 3}}
+{1\over \delta^{2}}\Big\|{\partial {\cal V}\over \partial s_\alpha}-\GR\Gt_\alpha\Big\|^2_{(L^2(\omega))^3}\ \le {C\over \delta^3}{\cal E}_\delta(\Gv).\cr}$$ }
\noindent{\Ggras Proof. } First of all there exists a positive constant $C$ independent  of $\delta$ such that
$$\left\{\eqalign{
&\delta\Big\|{\partial \GR\over \partial s_1} \Gn \cdot\GR\Gt_1\Big\|_{L^2(\omega)}+\delta\Big\| {\partial \GR\over \partial s_1} \Gn \cdot\GR\Gt_2 +{\partial \GR\over \partial s_2} \Gn \cdot\GR\Gt_1\Big\|_{L^2(\omega)}+\delta\Big\|{\partial \GR\over \partial s_2} \Gn \cdot\GR\Gt_2\Big\|_{L^2(\omega)}\cr
&+ ||{\cal Z}_{11}||_{L^2(\omega)}+ ||{\cal Z}_{12}||_{L^2(\omega)}+ ||{\cal Z}_{22}||_{L^2(\omega)}\le {C\over \delta^{1/2}}||\widehat{E}(\Gv)||_{(L^2(\Omega_\delta))^{3\times 3}},\cr
&\Big\|\GR^T{\partial\overline{v} \over \partial s_3} \cdot \Gt_\alpha+ {\cal Z}_{3\alpha}|\Big\|_{L^2(\Omega_\delta)}+\Big\|\GR^T{\partial\overline{v} \over \partial s_3} \cdot \Gn\Big\|_{L^2(\Omega_\delta)}\le C||\widehat{E}(\Gv)||_{(L^2(\Omega_\delta))^{3\times 3}}.\cr}\right.\leqno(5.7)$$ We use the definition of $\D_\delta$  to estimate the field $\GR^T\overline{v}\cdot\Gt_\alpha $. Introducing the function $\ds\GR^T\overline{v}\cdot\Gt_\alpha+s_3{\cal Z}_{3\alpha}$, using   PoincarŽ-Wirtinger's inequality and the first condition on $\GR^T\overline{v}$ in $\D_\delta$ give 
$$\Big\|\GR^T\overline{v}\cdot\Gt_\alpha+s_3{\cal Z}_{3\alpha}\Big\|_{L^2(\Omega_\delta)}+\Big\|\GR^T\overline{v}\cdot\Gn\Big\|_{L^2(\Omega_\delta)}\le C\delta||\widehat{E}(\Gv)||_{(L^2(\Omega_\delta))^{3\times 3}}.\leqno(5.8)$$ Now we use the second condition on $\overline{v}\cdot\Gt_\alpha$ (in the definition of $\D_\delta$) in the above estimates  and again (5.7) to get the estimates on $\ds\GR^T{\partial\overline{v}\over \partial s_3}$ and ${\cal Z}_{3\alpha}$
$$\sum_{\alpha=1}^2\Big\{\Big\|\GR^T{\partial\overline{v}\over \partial s_3}\cdot\Gt_\alpha\Big\|_{(L^2(\Omega_\delta))^3}+\delta^{1/2}||{\cal Z}_{3\alpha}||_{L^2(\omega)}\Big\}+\Big\|\GR^T{\partial\overline{v}\over \partial s_3}\cdot\Gn\Big\|_{(L^2(\Omega_\delta))^3}\le C ||\widehat{E}(\Gv)||_{(L^2(\Omega_\delta))^{3\times 3}}.$$ Finally (5.8) gives the $L^2$ estimate on $\overline{v}$. Let us notice that due to the last condition on $\overline{v}$ in $\D_\delta$, we obtain the same estimates that in the case where  $\Gv$ satisfies the limit kinematic condition $\ds {\partial{\cal V}\over \partial s_\alpha }=\GR\Gt_\alpha$.

\noindent There exist two antisymmetric matrices $\GA_1$ and $\GA_2$ in $(L^2(\omega))^{3\times 3}$  such that
$${\partial \GR\over \partial s_1}=\GR\GA_1\qquad\qquad {\partial \GR\over \partial s_2}=\GR\GA_2.$$ From (5.7) we get
$$\big\| \GA_1 \Gn \cdot\Gt_1\big\|_{L^2(\omega)}+\big\|\GA_1 \Gn \cdot\Gt_2 +\GA_2 \Gn \cdot\Gt_1\big\|_{L^2(\omega)}+\big\| \GA_2 \Gn \cdot\Gt_2\big\|_{L^2(\omega)}\le  {C\over \delta^{3/2}}||\widehat{E}(\Gv)||_{(L^2(\Omega_\delta))^{3\times 3}}.$$ Besides there exists a positive constant such 
$$\eqalign{
&||\GA_1||_{(L^2(\omega))^{3\times 3}}+||\GA_2||_{(L^2(\omega))^{3\times 3}}\cr
\le &C\Big\{\big\| \GA_1 \Gn \cdot\Gt_1\big\|_{L^2(\omega)}+\big\|\GA_1 \Gn \cdot\Gt_2 +\GA_2 \Gn \cdot\Gt_1\big\|_{L^2(\omega)}+\big\| \GA_2 \Gn \cdot\Gt_2\big\|_{L^2(\omega)}+||\GA_1\Gt_2-\GA_2\Gt_1||_{(L^2(\omega))^{3}}\Big\}.\cr}$$ Hence we get
$$\Big\|{\partial \GR\over \partial s_1}\Big\|_{(L^2(\omega))^{3\times 3}}+\Big\|{\partial \GR\over \partial s_2}\Big\|_{(L^2(\omega))^{3\times 3}}\le C\Big\{{1\over \delta^{3/2}}||\widehat{E}(\Gv)||_{(L^2(\Omega_\delta))^{3\times 3}}+ \Big\|{\partial \GR\over \partial s_1} \Gt_2-{\partial \GR\over \partial s_2} \Gt_1\Big\|_{(L^2(\omega))^{3}}\Big\}.$$   Due to the  estimates concerning the ${\cal Z}_{i\alpha}$ and the definition of ${\cal E}_\delta(\Gv)$ we finally obtain
$$\Big\|{\partial {\cal V}\over \partial s_\alpha}-\GR\Gt_\alpha\Big\|^2_{(L^2(\omega))^3}\le {C\over \delta}{\cal E}_\delta(\Gv).$$ \fin
We define now the set of the admissible triplets
$$\D_{\delta,\gamma_0}=\Big\{ \Gv=({\cal V} , \GR , \overline{v})\in\D_\delta \;|\;   {\cal V}=\phi,\enskip \GR=\GI_3\quad \hbox{on }\gamma_0\;\Big\}.$$  Notice that the triplet $\GI_d=(\phi,\GI_3,0)$ belongs to $\D_{\delta,\gamma_0}$ and it is associated to the deformation $v=I_d$. 

 In some sense, the following corollary gives two  Korn's type inequalities on the set $\D_{\delta,\gamma_0}$ with respect to  the quantity ${\cal E}_\delta(\Gv)$, the more accurate of which depending on the order of  ${\cal E}_\delta(\Gv)$.
\smallskip
\noindent {\Ggras Corollary 5.3. } {\it There exists a positive constant $C$ which does not depend on $\delta$ such that for all $\Gv\in \D_{\delta,\gamma_0}$
$$\eqalign{
&||{\cal V}-\phi||^2_{(H^1(\omega))^{3}}+\|\GR-\GI_3\|^2_{(H^1(\omega))^{3\times 3}}\le {C\over  \delta^3} {\cal E}_\delta(\Gv),\cr
&||{\cal V}-\phi||^2_{(H^1(\omega))^{3}}\le C\Big(1+ {1\over \delta}{\cal E}_\delta(\Gv)\Big).\cr}$$}
\noindent{\Ggras Proof. } Recall that $\GR=\GI_3$ and ${\cal V}=\phi$ on $\gamma_0$, then from Proposition  5.1 we obtain
$$||\GR-\GI_3||^2_{(H^1(\omega))^{3\times 3}}\le {C \over  \delta^3}{\cal E}_\delta(\Gv).$$ Using   the above estimate and again Proposition  5.1   we obtain the first estimate on ${\cal V}-\phi$ (recall that $\ds \Gt_\alpha={\partial \phi\over \partial s_\alpha}$). To obtain the second estimate on ${\cal V}-\phi$, notice that $||\GR-\GI_3||^2_{(L^2(\omega))^{3\times 3}}\le C $.\fin
\noindent{\Ggras 6. Elastic shells }
  \medskip
In this section  we consider a shell made of an elastic material. Its thickness $2\delta$ is fixed and belongs to $]0,2\delta_0]$. The local energy $W\ \ :\ \  \GS_3\longrightarrow \R^+$ is a  continuous function of symmetric matrices which satisfies the following assumptions which are similar to those adopted in [14], [15] and [16] (the reader is also referred to [8] for general introduction to elasticity)
$$\leqalignno{
&\exists c>0 \quad \hbox{such  that }\quad \forall E\in\GS_3\quad  W(E)\ge c|||E|||^2, &(6.1)\cr
&\forall\varepsilon >0,\quad \exists \theta >0,\quad \hbox{such  that }\quad \forall E\in\GS_3\quad |||E|||\le \theta\; \Longrightarrow \; |W(E)-Q(E)|\le \varepsilon |||E|||^2, &(6.2)\cr}$$
where $Q$ is a positive quadratic form defined on the set of $3\times 3$ symmetric matrices. Remark that $Q$ satisfies (6.1) with the same constant $c$. 

Still following [8], for any $3\times3$ matrix $F$, we set 
$$\widehat{W}(F)=\left\{\eqalign{
&W\Big({1\over 2}(F^TF-\GI_3)\Big) \quad\hbox{if} \quad \det(F)>0\cr
&+\infty \hskip 2.4cm \hbox{if} \quad \det(F)\le0.\cr}\right.\leqno(6.3)$$

Remark that due to (6.1), (6.3) and to the inequality $|||F^TF-\GI_3|||\ge dist(F,SO(3))$ if $\det (F)>0$, we have for any $3\times 3$ matrix $F$
$$\widehat{W}(F)\ge {c\over 4} dist (F,SO(3))^2. \leqno (6.4)$$
 
\smallskip
\noindent {\ggras Remark 6.1. } As a classical example of a local elastic energy satisfying the above assumptions, we mention the following St Venant-Kirchhoff's law (see [8]) for which
$$\widehat{W}(F)=\left\{\eqalign{
&{\lambda\over 8}\big(tr(F^TF-\GI_3) \big)^2+{\mu\over 4}tr\big((F^TF-\GI_3)^2\big)\quad\hbox{if}\quad \det(F)>0\cr
&+\infty\qquad \hbox{if}\qquad \det(F)\le 0.\cr}\right.$$ 

In order to take into account the boundary condition on the admissible deformations  we introduce the space 
$$\GU_\delta=\Bigl\{v\in (H^{1}({\cal Q}_\delta))^3\; |\; v=I_d\quad \hbox{on}\quad \Gamma_{0,\delta}\Big\}.\leqno(6.5)$$ 

Let $\kappa\ge1$. Now we assume that the shell is submitted to applied body forces $f_{\kappa,\delta}\in(L^2(\Omega_\delta))^3$ and we define   the total energy  
$J_{\kappa,\delta}(v)$\footnote{*}{ For later convenience, we have added the term $\ds \int_{{\cal Q}_\delta}f_{\kappa,\delta}(x)\cdot  I_d(x)dx$
to the usual standard energy, indeed this does not affect the  minimizing problem for  $J_{\kappa,\delta}$. }
 over $\GU_\delta$ by
$$J_{\kappa,\delta}(v)=\int_{{\cal Q}_\delta}\widehat{W}(\nabla_x v)(x)dx-\int_{{\cal Q}_\delta}f_{\kappa,\delta}(x)\cdot (v(x)-I_d(x))dx.\leqno(6.6)$$
To introduce the scaling on $f_{\kappa,\delta}$, let us consider   $f$ and $g$ in $(L^2(\omega))^3$ and  assume  that the force  $f_{\kappa,\delta}$ is given  by
$$f_{\kappa,\delta}(x)=
\delta^{\kappa^{'}} f(s_1,s_2)+\delta^{\kappa^{'}-2}s_3 g (s_1,s_2)
\qquad \hbox{for a.e.}\enskip x=\Phi(s)\in {\cal Q}_\delta.\leqno(6.7)$$ where
$$\kappa^{'}=\left\{\eqalign{
&2\kappa-2\quad \hbox{if }\enskip 1\le\kappa\le 2,\cr
&\kappa\hskip 11mm \hbox{if }\enskip \kappa\ge 2.\cr}\right.\leqno(6.8)$$ 

Notice that $J_{\kappa,\delta}(I_d)=0$. So,  in order to minimize $J_{\kappa,\delta}$ we only need to  consider  deformations $v $ of $\GU_{\delta}$ such that  
$J_{\kappa,\delta}(v )\le 0$. 

Now from (6.1), (6.3), (6.4), the two Korn's type inequalities (3.1)-(3.2), the assumption (6.7) of the body forces   and the definition (6.8) of $\kappa^{'}$, we obtain the following bound   for $||\hbox{dist}(\nabla_x v, SO(3))||_{L^2({\cal Q}_\delta)}$
$$||\hbox{dist}(\nabla_x v , SO(3))||_{L^2({\cal Q}_\delta)}\le C\delta^{\kappa-1/2}\quad \hbox{and}\quad \int_{{\cal Q}_\delta}f_{\kappa,\delta}\cdot(v -I_d)\le C\delta^{2\kappa-1} \leqno(6.9)$$ which in turn imply that
$$c\delta^{2\kappa-1}\le J_{\kappa,\delta}(v)\le 0.\leqno(6.10) $$ 

Again from  (6.3)-(6.4) and  the estimates (6.9)  we deduce
$$\eqalign{
{c\over 4}||(\nabla_x v )^T\nabla_x v -\GI_3||^2_{(L^2({\cal Q}_\delta))^{3\times 3}} & \le J_{\kappa,\delta}(v)+\int_{{\cal Q}_\delta}f_{\kappa,\delta}\cdot(v -I_d)\le C\delta^{2\kappa-1}.\cr}$$ Hence, the following estimate of the Green-St Venant's tensor:
$$\big\|{1\over 2}\big\{(\nabla_x v )^T\nabla_x v -\GI_3\big\}\big\|_{(L^2({\cal Q}_\delta))^{3\times 3}}\le  C\delta^{\kappa-1/2}.$$ We deduce from the above inequality that $v\in (W^{1,4}({\cal Q}_\delta))^3$  with
$$||\nabla_x v ||_{(L^4({\cal Q}_\delta))^{3\times 3}}\le  C\delta^{1\over 4}.$$

We set
$$m_{\kappa,\delta}=\inf_{v\in \GU_\delta}J_{\kappa,\delta}(v).$$
As a consequence of  (6.10) we  have 
$$c\le {m_{\kappa,\delta}\over \delta^{2\kappa-1}}\le 0.$$ 

 In general, a minimizer of $J_{\kappa,\delta}$ does not exist  on $\GU_\delta$. In what follows, we replace the elastic functional  $\ds v\longmapsto J_{\kappa,\delta}(v)$ on $\GU_\delta$  by a simplified functional defined on $\D_\delta$  which admits a minimum.
 
\medskip
{\it From now on we assume $\kappa>1$.}
\medskip
\noindent{\Ggras 7. The simplified  elastic model for shells }
 \medskip

The aim of this section is to define a functional $J^s_{\kappa,\delta}$ on the set $\D_{\delta,\gamma_0}$, which will appear as a simplification of the total energy $J_{\kappa,\delta}$ defined on the set $\GU_\delta$. In order to perform this task, we use the results of Section 5 and we proceed in three steps. Let us first consider an admissible deformation $v$   satisfying (3.3), decomposed as in (3.4) and such that  $J_{\kappa,\delta}(v)\le 0$.  It is convenient to  express the energy  $J_{\kappa,\delta}(v)$ over the  domain $\Omega_\delta$ 
$$\eqalign{
J_{\kappa,\delta}(v)= \int_{\Omega_\delta}{W}\Big({1\over 2 }\big((\nabla_x v )^T\nabla_x v-\GI_3\big)\Big)\det \big(\Gt_1 +s_3{\partial \Gn\over \partial s_1} | \Gt_2Ê+s_3{\partial \Gn\over \partial s_2} | \Gn\big)ds_1ds_2ds_3\cr
-\int_{\Omega_\delta} \big(\delta^{\kappa^{'}}f +\delta^{\kappa^{'}-2}s_3 g\big)\cdot \big(v-I_d\big)\det \big(\Gt_1 +s_3{\partial \Gn\over \partial s_1} | \Gt_2Ê+s_3{\partial \Gn\over \partial s_2} | \Gn\big)ds_1ds_2ds_3.\cr}\leqno(7.1)$$

The triplet associated to $v$ by the decomposition (3.4)  is denoted $\Gv=\big({\cal V}, \GR,\overline{v}\big)$.
The following estimate has been proved in Section 6
$$\Big\|{1\over 2} \big\{(\nabla_x v )^T\nabla_x v -\GI_3\big\}\Big\|_{(L^2(\Omega_\delta))^{3\times 3}}\le  C\delta^{\kappa-1/2}.$$
Then, for all $\theta>0$,  the set $A_\delta^\theta=\{s\in\Omega; |||\Pi _\delta \big((\nabla_xv_\delta)^T\nabla_x v_\delta-\GI_3\big)(s)||| \geq \theta\}$ has a measure satisfying
$$ {\rm meas}(A_\delta^\theta)\le C{\delta^{2\kappa-2}\over \theta^2}.$$ Now, according to assumptions (6.2) and $\kappa>1$ and the above estimate, in the first term of the total energy  $J_{\kappa,\delta}(v)$ we replace the quantity $\ds{W}\Big({1\over 2 } \big((\nabla_x v )^T\nabla_x v-\GI_3\big)\Big)$ by $\ds{Q}\Big({1\over 2 } \big((\nabla_x v )^T\nabla_x v-\GI_3\big)\Big)$. Following the analysis of Section 5, we then replace $\ds{Q}\Big({1\over 2 } \big((\nabla_x v )^T\nabla_x v-\GI_3\big)\Big)$ by $Q\big((\Gt_1\,|\, \Gt_2\,|\, \Gn)^{-T} \widehat{E}(\Gv)(\Gt_1\,|\, \Gt_2\,|\, \Gn)^{-1}\big)$ where $ \widehat{E}(\Gv)$ is defined by (5.3) and (5.5). At last, we replace $\ds \det \big(\Gt_1 +s_3{\partial \Gn\over \partial s_1} | \Gt_2Ê+s_3{\partial \Gn\over \partial s_2} | \Gn\big)$ by $\det \big(\Gt_1  | \Gt_2  | \Gn\big)$.
Setting for all $3\times 3$ symmetric matrix $F$
$$W^s(F)=Q\Big((\Gt_1\,|\, \Gt_2\,|\, \Gn)^{-T}\, F\, (\Gt_1\,|\, \Gt_2\,|\, \Gn)^{-1}\Big)\leqno(7.2)$$ all the above considerations lead us to replace the first term in the right hand side of (7.1) by
$$\ds \int_{\Omega_\delta}W^s\big( \widehat{E}(\Gv)\big)\det (\Gt_1 | \Gt_2Ê| \Gn)ds_1ds_2ds_3.\leqno(7.3)$$ 

\noindent Observe now the term involving the forces in (7.1). We have
$$ \eqalign{
\Big|&\int_{\Omega_\delta} \big(\delta^{\kappa^{'}}f +\delta^{\kappa^{'}-2}s_3 g\big)\cdot  \big(v-I_d\big)\det \big(\Gt_1 +s_3{\partial \Gn\over \partial s_1} | \Gt_2Ê+s_3{\partial \Gn\over \partial s_2} | \Gn\big)ds_1ds_2ds_3\cr
&-2\delta^{\kappa^{'}+1}\int_\omega \Big[f  \cdot({\cal V}-\phi) +{1\over 3}g  \cdot\big(\GR -\GI_3\big)\Gn\Big]\det (\Gt_1 | \Gt_2Ê| \Gn)ds_1 ds_2\cr
&- {2\over 3}\delta^{\kappa^{'}+1}\int_\omega g\cdot({\cal V}-\phi) \Big[\det\Big({\partial \Gn \over \partial s_1} |\Gt_2\; |\; \Gn\Big)+\det\Big(\Gt_1 | {\partial \Gn \over \partial s_2} \; |\; \Gn\Big)\Big]ds_1ds_2\Big|\cr
\le &C\delta^{\kappa^{'}+2}\big(||f||_{(L^2(\omega))^3}+||g||_{(L^2(\omega))^3}\big)\big(||{\cal V}-\phi||_{(L^2(\omega))^3}+||\GR-\GI_3||_{(L^2(\omega))^{3\times 3}}+{1\over \delta^{5/2}}||\overline{v}||_{(L^2(\Omega_\delta))^3}\big).\cr}$$  
Then, in view of the first estimate in (3.5)  we replace the term involving the forces by  
$${\cal L}_{\kappa,\delta}({\cal V}, \GR)=\delta^{\kappa^{'}+1}{\cal L}({\cal V}, \GR)\leqno(7.4)$$
where
$$\eqalign{
{\cal L}({\cal V}, \GR) &=2\int_\omega \Big[f \cdot({\cal V}-\phi)+{1\over 3}g \cdot(\GR-\GI_3)\Gn\Big]\det (\Gt_1 | \Gt_2Ê| \Gn)ds_1 ds_2\cr
&+ {2\over 3}\int_\omega g \cdot({\cal V}-\phi)\Big[\det\Big({\partial \Gn \over \partial s_1} |\Gt_2\; |\; \Gn\Big)+\det\Big(\Gt_1 | {\partial \Gn \over \partial s_2} \; |\; \Gn\Big)\Big]ds_1ds_2.\cr}$$ 
At the end of this first step, we obtain a simplified energy for a deformation $v\in\GU_\delta$ which satisfies (3.3) and $J_{\kappa,\delta}(v)\le 0$
$$J^{simpl}_{\kappa,\delta}(v)=\int_{\Omega_\delta}W^s\big(\widehat{E}(\Gv)\big)\det (\Gt_1 | \Gt_2Ê| \Gn)ds_1ds_2ds_3 - \delta^{\kappa^{'}+1}{\cal L}({\cal V}, \GR).$$
Indeed the energy $J^{simpl}_{\kappa,\delta}(v)$ can be seen as a functional of $\Gv$ defined over $ \D_{\delta,\gamma_0}$ since we have already notice that $\widehat{E}(\Gv)$ belongs to $(L^2(\Omega_\delta))^{3\times 3}$. As a consequence, in a second step we are in a position to extend the above energy to the whole set $\D_{\delta,\gamma_0}$ and to put
$$\forall \Gv\in \D_{\delta,\gamma_0},\qquad  J^{\bf simpl}_{\kappa,\delta}(\Gv)=\int_{\Omega_\delta}W^s\big(\widehat{E}(\Gv)\big)\det (\Gt_1 | \Gt_2Ê| \Gn)ds_1ds_2ds_3 - \delta^{\kappa^{'}+1}{\cal L}({\cal V}, \GR).$$
As observed in Section 5, the functional $J^{\bf simpl}_{\kappa,\delta}$ is not coercive on $\D_{\delta,\gamma_0}$. In a third step, in view of Proposition 5.2 and in order to obtain the coerciveness of the simplified energy, the two terms $\ds\delta^3\Big\|{\partial \GR\over \partial s_1} \Gt_2-{\partial \GR\over \partial s_2} \Gt_1\Big\|^2_{(L^2(\omega))^{3}}$, 
$\ds\delta\Big\|{\partial{\cal V}\over \partial s_1 }\cdot\GR\Gt_2- {\partial{\cal V}\over \partial s_2}\cdot\GR\Gt_1\Big\|^2_{L^2(\omega)}$ are added to  $J^{\bf simpl}_{\kappa,\delta}$.
\smallskip

Using all the above considerations, we are able to define the simplified elastic energy on $\D_{\delta,\gamma_0}$ by setting  for any  $\Gv$ in $\D_{\delta,\gamma_0}$ 
$$\left\{\eqalign{
J^s_{\kappa,\delta}(\Gv)= &\int_{\Omega_\delta}W^s\big(\widehat{E}(\Gv)\big)\det (\Gt_1 | \Gt_2Ê| \Gn)ds_1ds_2ds_3+ \delta^3\Big\|{\partial \GR\over \partial s_1} \Gt_2-{\partial \GR\over \partial s_2} \Gt_1\Big\|^2_{(L^2(\omega))^{3}}\cr
&+\delta\Big\|{\partial{\cal V}\over \partial s_1 }\cdot\GR\Gt_2- {\partial{\cal V}\over \partial s_2}\cdot\GR\Gt_1\Big\|^2_{L^2(\omega)}  - \delta^{\kappa^{'}+1}{\cal L}({\cal V}, \GR).\cr}\right.\leqno(7.5)$$

\smallskip
The end  of this section is dedicated to show that the functional  $J^s_{\kappa,\delta}$ admits a minimizer on $\D_{\delta,\gamma_0}$.  Let $\Gv$ be in $\D_{\delta,\gamma_0}$  we have
$$\Big|{\cal L}({\cal V}, \GR)\Big| \le C \big(||f||_{(L^2(\omega))^3}+||g||_{(L^2(\omega))^3}\big)\big(||{\cal V}-\phi||_{(L^2(\omega))^3}+||\GR-\GI_3||_{(L^2(\omega))^{3\times 3}}\big). \leqno(7.6)$$ 
The quadratic form $Q$ being positive, the definition (5.6) of ${\cal E}_\delta(\Gv)$ and (7.5)-(7.6) give
$$C {\cal E}_\delta(\Gv)-C\delta^{\kappa^{'}+1}\big(||f||_{(L^2(\omega))^3}+||g||_{(L^2(\omega))^3}\big)\big(||{\cal V}-\phi||_{(L^2(\omega))^3}+||\GR-\GI_3||_{(L^2(\omega))^{3\times 3}}\big)\le J^s_{\kappa,\delta}(\Gv).$$ Now  thanks to Corollary 5.3 and (6.8), we get, if $J^s_{\kappa,\delta}(\Gv)\le 0\big(=J^s_{\kappa,\delta}(\GI_d)\big)$ 
$${\cal E}_\delta(\Gv)\le C\delta^{2\kappa-1}(||f||_{(L^2(\omega))^3}+||g ||_{(L^2(\Omega))^3})^2.\leqno(7.7)$$ Hence, there exists  a  constant $c$  which does not depend  on $\delta$  such that for any $\Gv\in \D_{\delta,\gamma_0}$ satisfying $ J^s_{\kappa,\delta}(\Gv) \le 0$,  we have
$$c\delta^{2\kappa-1}\le J^s_{\kappa,\delta}(\Gv) .$$  We set
$$m^s_{\kappa,\delta}=\inf_{\Gv\in \D_{\delta,\gamma_0}} J^s_{\kappa,\delta}(\Gv).\leqno(7.8)$$ As a consequence of the  above inequality, we  have 
$$c\le {m^s_{\kappa,\delta}\over \delta^{2\kappa-1}}\le 0.$$ In the following theorem we prove that for $\kappa$ and  $\delta$ fixed the minimization problem (7.8) has at least a  solution.
\smallskip
\noindent{\Ggras Theorem 7.1. }{\it There exists $\Gv_\delta\in \D_{\delta,\gamma_0}$ such that
$$m^s_{\kappa,\delta}=J^s_{\kappa,\delta}(\Gv_\delta) =\min_{v\in \D_{\delta,\gamma_0}}J^s_{\kappa,\delta}(\Gv).\leqno(7.9)$$}
\noindent{\ggras  Proof. } Since  $J^s_{\kappa,\delta}(\GI_d) = 0$, we can consider a minimizing sequence  $\Gv_n$ in $\D_{\delta,\gamma_0}$ such that $J^s_{\kappa,\delta}(\Gv_n) \le 0$ and
$$m^s_{\kappa,\delta}=\lim_{n\to+\infty}J^s_{\kappa,\delta}(\Gv_n).$$ From (7.7) we get
$$ {\cal E}_\delta(\Gv_n)\le C\delta^{2\kappa-1}(||f_{(L^2(\omega))^3}+||g ||_{(L^2(\Omega))^3})^4.$$ Thanks to Corollary 5.3 and Proposition 5.2, the above estimate show that there exists a  subsequence still denoted $n$ such that (recall that $||\GR_n||_{(L^\infty(\omega))^{3\times 3}}=\sqrt 3$)
$$\eqalign{
{\cal V}_n & \rightharpoonup {\cal V}_\delta\quad \hbox{weakly in } (H^1(\omega))^3\cr
\GR_n & \rightharpoonup \GR_\delta\quad \hbox{weakly in } (H^1(\omega))^{3\times 3}\cr
\GR_n & \longrightarrow \GR_\delta\quad \hbox{strongly in } (L^2(\omega))^{3\times 3}\enskip\hbox{and a.e. in } \omega\cr
\overline{v}_n&\rightharpoonup \overline{v}_\delta\quad \hbox{weakly in } (L^2(\omega ; H^1(-\delta,\delta)))^3.\cr}$$ Then setting $\Gv_\delta=({\cal V}_\delta,\GR_\delta,\overline{v}_\delta)\in \D_{\delta,\gamma_0}$,  we get
$$\eqalign{
&\widehat{E}(\Gv_n) \rightharpoonup \widehat{E}(\Gv_\delta)\quad \hbox{weakly in } (L^2(\Omega_\delta))^{3\times 3},\cr
&{\partial{\cal V}_n\over \partial s_1 }\cdot\GR_n\Gt_2- {\partial{\cal V}_n\over \partial s_2}\cdot\GR_n\Gt_1\rightharpoonup{\partial{\cal V}_\delta\over \partial s_1 }\cdot\GR_\delta\Gt_2- {\partial{\cal V}_\delta\over \partial s_2}\cdot\GR_\delta\Gt_1\quad \hbox{weakly in } (L^2(\omega))^{3}.\cr}$$ Now, passing to the limit inf  in $J^s_{\kappa,\delta}(\Gv_n)$, we obtain
$$m^s_{\kappa,\delta}\le J^s_{\kappa,\delta}(\Gv_\delta) \le \liminf_{n\to+\infty}J^s_{\kappa,\delta}(\Gv_n)= \lim_{n\to+\infty}J^s_{\kappa,\delta}(\Gv_n)= m^s_{\kappa,\delta}.$$\fin
\smallskip
\noindent {\Ggras  8.  Asymptotic behavior of the simplified model.  Case $\kappa=2$.}
\smallskip
 In this section  we study the asymptotic behavior of the sequence $(\Gv_\delta)$ of minimizer given in Theorem 7.1 and we characterize  the limit of the minima $\ds {m^s_{2,\delta}\over \delta^3}$ as a minimum of a new functional. AS usual, to perform this task, we work on the fixed domain $\Omega$  and we use the operator $\Pi_\delta$ defined in Section 4. We denote $\D$ the following closed subset of $\D_{1,\gamma_0}$ (i.e. $\D_{\delta,\gamma_0}$ for $\delta=1$ or  $\D_{1,\gamma_0}=\Pi_\delta \big( \D_{\delta,\gamma_0}\big)$):
$$\D=\Big\{ \Gv=({\cal V} , \GR , \overline{V})\in    \D_{1,\gamma_0}\;\;|\; {\partial {\cal V}\over \partial s_\alpha}=\GR\Gt_\alpha \Big\}.$$ 
Notice that ${\cal V}\in (H^2(\omega))^3$.  Then we define the following functional over $\D$
 $${\cal J}_2(\Gv)= \int_{\Omega} Q\Big((\Gt_1\,|\, \Gt_2\,|\, \Gn)^{-T}\GE(\Gv) (\Gt_1\,|\, \Gt_2\,|\, \Gn)^{-1}\Big)\det (\Gt_1 | \Gt_2Ê| \Gn)-{\cal L}({\cal V},\GR).\leqno(8.1)$$
where
$$\GE(\Gv)=\pmatrix{
\displaystyle  S_3{\partial \GR\over \partial s_1} \Gn \cdot\GR\Gt_1 & \displaystyle  S_3 {\partial \GR\over \partial s_1} \Gn \cdot\GR\Gt_2  &\displaystyle  {1\over 2}\GR^T{\partial\overline{V} \over \partial S_3} \cdot \Gt_1 \cr
 * & \displaystyle  S_3{\partial \GR\over \partial s_2} \Gn \cdot\GR\Gt_2  &\displaystyle  {1\over 2}\GR^T{\partial\overline{V} \over \partial S_3} \cdot\Gt_2\cr
 * & *&  \displaystyle \GR^T{\partial\overline{V} \over \partial S_3} \cdot\Gn \cr}\leqno(8.2)$$   As in Theorem 7.1 we easily prove that there exists $\Gv_2=({\cal V}_2,\GR_2,\overline{V}_2)\in\D$ such that
$$m^s_2={\cal J}_2(\Gv_2)=\min_{\Gv\in \D}{\cal J}_2(\Gv).\leqno(8.3)$$ 
\noindent{\ggras  Theorem  8.1. }{\it We have 
$$m^s_2=\lim_{\delta\to 0}{m^s_{2,\delta}\over \delta^{3}}.$$
Moreover, let $\Gv_\delta=({\cal V}_\delta,\GR_\delta,\overline{v}_\delta)\in \D_{\delta,\gamma_0}$ be a minimizer of the functional $J^s_{2,\delta}(\cdot )$, there exists a subsequence  still denoted $\delta$ such that
$$\left\{\eqalign{
{\cal V}_\delta& \longrightarrow {\cal V}_0\quad \hbox{strongly in } (H^1(\omega))^3,\cr
\GR_\delta & \longrightarrow  \GR_0\quad \hbox{strongly in } (H^1(\omega))^{3\times 3},\cr
{1\over \delta}  {\cal Z}_{i\beta,\delta}&\longrightarrow 0\quad \hbox{strongly in } L^2(\omega),\cr
{1\over \delta^2}\Pi_\delta(\overline{v}_\delta)&\longrightarrow \overline{V}_0\quad \hbox{strongly in } (L^2(\omega ; H^1(-1,1)))^3.\cr}\right.\leqno(8.4)$$ The triplet $\Gv_0=({\cal V}_0,\GR_0,\overline{V}_0)$ belongs to $\D$ and we have
$$m^s_2={\cal J}_{2}(\Gv_0).$$ }
 \noindent{\ggras  Proof. } For all $\Gv=({\cal V},\GR, \overline{V})\in \D$, we have $({\cal V},\GR,  \overline{v}_\delta)\in\D_{\delta,\gamma_0}$ where
 $$ \overline{v}_\delta(s_1,s_2,s_3)=\delta^2\overline{V}(s_1,s_2,{s_3\over \delta}\Big)\qquad \hbox{for a.e. } (s_1,s_2,s_3)\in \Omega_\delta.$$ Using the fact that $\Gv \in \D$, which implies that  $\ds {\partial \GR\over \partial s_1} \Gn \cdot\GR\Gt_2={\partial \GR\over \partial s_2} \Gn \cdot\GR\Gt_1$, we have
 $${J_{2,\delta}^s\big({\cal V},\GR,  \overline{v}_\delta\big)\over \delta^3}= \int_{\Omega} Q\Big((\Gt_1\,|\, \Gt_2\,|\, \Gn)^{-T}\GE(\Gv) (\Gt_1\,|\, \Gt_2\,|\, \Gn)^{-1}\Big)\det (\Gt_1 | \Gt_2Ê| \Gn)-{\cal L}({\cal V},\GR)={\cal J}_2(\Gv).\leqno(8.5)$$ Then, taking the minimum in the right hand side w.r.t. $\Gv \in\D$, we immediately deduce that  $\ds {m^s_{2,\delta}\over \delta^{3}}\le m^s_2$. 
\smallskip
\noindent We recall that   $\Gv_\delta=\big({\cal V}_\delta, \GR_\delta, \overline{v}_\delta\big)\in D_{\delta,\gamma_0}$ is a minimizer of $J^s_\delta$ 
$$c\le{m^s_{2,\delta}\over \delta^{3}}={J^s_{2,\delta}(\Gv_\delta)\over \delta^{3}}=\min_{\Gv\in \D_{\delta,\gamma_0}}{J^s_{2,\delta}(\Gv)\over \delta^{3}} $$ and moreover with (7.7)
$${\cal E}_\delta(\Gv_\delta) \le C\delta^3(||f_{(L^2(\omega))^3}+||g ||_{(L^2(\Omega))^3})^2.$$Thanks to the estimates in Proposition 5.2, Corollary 5.3 and the above estimate we can extract a subsequence  still denoted $\delta$ such that
$$\left\{\eqalign{
{\cal V}_\delta& \longrightarrow {\cal V}_0\quad \hbox{strongly in } (H^1(\omega))^3,\cr
\GR_\delta & \rightharpoonup \GR_0\quad \hbox{weakly in } (H^1(\omega))^{3\times 3}\quad \hbox{and a.e. in }  \omega,\cr
{1\over \delta^2}\Pi_\delta(\overline{v}_\delta)&\rightharpoonup \overline{V}_0\quad \hbox{weakly in } (L^2(\omega ; H^1(-1,1)))^3,\cr
{1\over \delta}  {\cal Z}_{i\alpha,\delta}& \rightharpoonup {\cal Z}_{i\alpha,0}\quad \hbox{weakly in } L^2(\omega),\cr
  \Big({\partial \GR_\delta\over \partial s_1}\Gt_2&-{\partial \GR_\delta\over \partial s_2}\Gt_1\Big) \rightharpoonup Y\quad \hbox{weakly in }( L^2(\omega))^{3},\cr
{1\over \delta} \Big({\partial {\cal V}_\delta\over \partial s_1}\cdot\GR_\delta\Gt_2&-{\partial {\cal V}_\delta\over \partial s_2}\cdot \GR_\delta\Gt_1\Big) \rightharpoonup X\quad \hbox{weakly in } L^2(\omega).\cr}\right.\leqno(8.6)$$ 
Then  from the fifth convergence  we obtain  
$\ds {\partial {\cal V}_0\over \partial s_\alpha}=\GR_0\Gt_\alpha$. So we have ${\cal V}_0\in (H^2(\omega))^3$ and $\Gv_0=({\cal V}_0,\GR_0,\overline{V}_0)$ belongs to $\D$.  From the above convergences, and upon extracting another subsequence, we also get
$${1\over \delta}\Pi_\delta\big(\widehat{E}(\Gv_\delta)\big)\rightharpoonup\GE_0 \quad\hbox{weakly in }( L^2(\Omega))^{3\times 3}$$ where
 $$ \GE_0=\pmatrix{
\displaystyle  S_3{\partial \GR_0\over \partial s_1} \Gn \cdot\GR_0\Gt_1+{\cal Z}_{11,0} & \displaystyle  S_3 {\partial \GR_0\over \partial s_1} \Gn \cdot\GR_0\Gt_2+ {\cal Z}_{12,0} &\displaystyle  {1\over 2}\GR^T_0{\partial\overline{W}_0 \over \partial S_3} \cdot \Gt_1 \cr
 * & \displaystyle  S_3{\partial \GR_0\over \partial s_2} \Gn \cdot\GR_0\Gt_2+{\cal Z}_{22,0}   &\displaystyle  {1\over 2}\GR^T_0{\partial\overline{W}_0 \over \partial S_3} \cdot\Gt_2\cr
 * & *&  \displaystyle  \GR^T_0{\partial\overline{W}_0\over \partial S_3} \cdot\Gn \cr}$$ 
with
 $$\overline{W}_0=\overline{V}_0+S_3{\cal Z}_{31,0}\GR_0\Gt^{'}_1+S_3{\cal Z}_{32,0}\GR_0\Gt^{'}_2.$$ Due to the expression of $J^s_{2,\delta}$ we have
$$\eqalign{
{J^s_{2,\delta}(\Gv_\delta)\over \delta^3}= &\int_{\Omega}Q\Big((\Gt_1\,|\, \Gt_2\,|\, \Gn)^{-T}\,{1\over \delta}\Pi_\delta(\widehat{E}(\Gv_\delta))\,(\Gt_1\,|\, \Gt_2\,|\, \Gn)^{-1}\Big)\det (\Gt_1 | \Gt_2Ê| \Gn)ds_1ds_2dS_3+ \Big\|{\partial \GR_\delta\over \partial s_1} \Gt_2-{\partial \GR_\delta\over \partial s_2} \Gt_1\Big\|^2_{(L^2(\omega))^{3}}\cr
&+{1\over \delta^2}\Big\|{\partial{\cal V}_\delta\over \partial s_1 }\cdot\GR_\delta\Gt_2- {\partial{\cal V}_\delta\over \partial s_2}\cdot\GR_\delta\Gt_1\Big\|^2_{L^2(\omega)}-{\cal L}({\cal V}_\delta, \GR_\delta).\cr}$$ 
With the convergences (8.6), since $Q$ is quadratic and thanks to the expression of ${\cal L}$,  we are in a position to pass to the limit-inf in the above equality which gives
$$\eqalign{
 \int_\Omega Q\Big((\Gt_1\,|\, \Gt_2\,|\, \Gn)^{-T}\GE_0 (\Gt_1\,|\, \Gt_2\,|\, \Gn)^{-1}\Big)\det (\Gt_1 | \Gt_2Ê| \Gn)ds_1ds_2dS_3+||X||^2_{L^2(\omega)}+||Y||^2_{(L^2(\omega))^3}-{\cal L}({\cal V}_0,\GR_0)\cr
 \le \liminf_{\delta\to 0}{J^s_{2,\delta}(\Gv_\delta)\over \delta^{3}}=\liminf_{\delta\to 0}{m^s_{2,\delta}\over \delta^{3}}.\cr}$$ Hence we get
$$\int_\Omega Q\Big((\Gt_1\,|\, \Gt_2\,|\, \Gn)^{-T}\GE_0 (\Gt_1\,|\, \Gt_2\,|\, \Gn)^{-1}\Big)\det (\Gt_1 | \Gt_2Ê| \Gn)ds_1ds_2dS_3-{\cal L}({\cal V}_0,\GR_0)\le \liminf_{\delta\to 0}{m^s_{2,\delta}\over \delta^{3}}.$$ 
First, notice that if $\Gv=({\cal V}, \GR, \overline{V})\in \D$ then $\overline{V}$ satisfies
$$ \int_{-1}^1{\partial \overline{V}\over \partial S_3}(s_1,s_2,S_3)\cdot\Gt_\alpha(s_1,s_2)(S^2_3-1)dS_3=0\qquad \hbox{for a.e. } (s_1,s_2)\in \omega.$$ Now we apply  Lemma A with $\ds\Ga=\Big({\partial \GR_0\over \partial s_1} \Gn \cdot\GR_0\Gt_1, {\partial \GR_0\over \partial s_1} \Gn \cdot\GR_0\Gt_2, {\partial \GR_0\over \partial s_2} \Gn \cdot\GR_0\Gt_2\Big)$, $\Gb=({\cal Z}_{11,0},  {\cal Z}_{12,0}, {\cal Z}_{22,0})$,  $\ds\Gc= \Big({1\over 2}\GR^T_0{\partial\overline{W}_0 \over \partial S_3} \cdot \Gt_1,{1\over 2}\GR^T_0{\partial\overline{W}_0 \over \partial S_3} \cdot \Gt_2, \GR^T_0{\partial\overline{W}_0 \over \partial S_3} \cdot \Gn\Big)$ and
with the  quadratic form defined by
$${\cal Q}_m(\Ga,\Gb,\Gc)=\int_{-1}^1 Q\Big((\Gt_1\,|\, \Gt_2\,|\, \Gn)^{-T}\GE_0 (\Gt_1\,|\, \Gt_2\,|\, \Gn)^{-1}\Big)\det (\Gt_1 | \Gt_2Ê| \Gn)dS_3\qquad \hbox{for a.e. } (s_1,s_2)\in\omega.$$ We obtain
$$\min_{\Gv\in \D}{\cal J}_2(\Gv)\le  \int_\Omega Q\Big((\Gt_1\,|\, \Gt_2\,|\, \Gn)^{-T}\GE_0 (\Gt_1\,|\, \Gt_2\,|\, \Gn)^{-1}\Big)\det (\Gt_1 | \Gt_2Ê| \Gn)ds_1ds_2dS_3-{\cal L}({\cal V}_0,\GR_0).\leqno(8.7)$$ Hence $\ds m^s_2\le \liminf_{\delta\to 0}{m^s_{2,\delta}\over \delta^{3}}$. Recall that we have $\ds {m^s_{2,\delta}\over \delta^{3}}\le m^s_2$, so we get
$$\lim_{\delta\to 0}{m^s_{2,\delta}\over \delta^{3}}=m^s_2.$$
Finally, from convergences (8.6)  we obtain  ${\cal Z}_{i\alpha,0}=0$, $X=Y=0$ and moreover we have the strong convergences in (8.4). \fin \smallskip
\noindent {\Ggras  9.   Justification of the simplified model. Case $\kappa=2$}
\medskip
In this section, the introduction of the simplified energy is justified in the sense that we prove that both the minima of the elastic energy and of the simplified energy have the same limit as $\delta$ tends to $0$.

\noindent{\ggras  Theorem  9.1. }{\it We have
$$\lim_{\delta\to 0}{m_{2,\delta}\over \delta^{3}}=\lim_{\delta\to 0}{m^s_{2,\delta}\over \delta^{3}}=m^s_2.$$}
\noindent{\ggras  Proof. }

\noindent{\ggras Step 1.} In this step we prove that $\ds m^s_2\le \liminf_{\delta\to 0}{m_{2,\delta}\over \delta^{3}}$. Let   $\big(v_\delta\big)_{0<\delta\le \delta_0 }$  be  a minimizing sequence of deformations belonging to $\GU_\delta$ and such that
$$\liminf_{\delta\to 0}{m_{2,\delta}\over \delta^{3}}=\lim_{\delta\to 0}{J_{2,\delta}(v_\delta)\over \delta^3}.\leqno(9.1)$$  From the estimates of Section 6 we get
$$\left\{\eqalign{
&||\hbox{dist}(\nabla_x v_\delta, SO(3))||_{L^2({\cal Q}_\delta)}\le C\delta^{3/2},\cr
& \big\|{1\over 2}\big\{\nabla_x v_\delta^T\nabla_x v_\delta-\GI_3\big\}\big\|_{(L^2({\cal Q}_\delta))^{3\times 3}}\le  C\delta^{3/2},\cr
 &||\nabla_x v_\delta||_{(L^4({\cal Q}_\delta))^{3\times 3}}\le  C\delta^{1/ 4}.\cr}\right.\leqno(9.2)$$ 
  We still denote by $\ds{\cal V}_\delta(s_1,s_2)={1\over 2\delta}\int_{-\delta}^\delta v_\delta(s_1,s_2,s_3)ds_3$ the mean of $v_\delta$ over the fibers of the shell. Upon extracting a subsequence (still indexed by $\delta$), the results of  [6] show   that there exist ${\cal V}\in (H^2(\omega))^3$, $\GR\in (H^1(\omega))^{3\times 3}$ with $\GR(s_1,s_2) \in SO(3)$ for a.e.  $(s_1,s_2) \in \omega$, ${\cal Z}_{\alpha\beta}\in L^2(\omega)$ and  $\overline{V} \in \big(L^2(\omega; H^1(-1,1)\big)^3 $ satisfying
  $$\int_{-1}^1\overline{V}(s_1,s_2,S_3)dS_3=0\qquad \hbox{ for a.e.  }\enskip (s_1,s_2) \in \omega,\qquad\qquad {\partial{\cal V} \over \partial s_\alpha}=\GR  \Gt _\alpha \leqno (9.3)$$ together with the boundaries conditions
${\cal V} =\phi$, $\GR=\GI_3$ on $\gamma_0,$ 
and with the following convergences
$$\left\{\eqalign{
&\Pi _\delta(v_\delta) \longrightarrow    {\cal V}  \quad \hbox{strongly in}\quad \big(H^1(\Omega)\big)^3,\cr
&\Pi _\delta (\nabla_x v_\delta)  \longrightarrow  \GR \quad \hbox{strongly in}\quad \big(L^2(\Omega)\big)^{3\times 3},\cr
&{\Pi_\delta(v_\delta-{\cal V}_\delta)\over \delta}\longrightarrow S_3(\GR-\GI_3)\Gn\quad \hbox{strongly in}\quad \big(L^2(\Omega)\big)^3,\cr
&{1\over 2\delta}\Pi _\delta \big((\nabla_xv_\delta)^T\nabla_x v_\delta-\GI_3\big) \rightharpoonup   (\Gt_1\,|\, \Gt_2\,|\, \Gn)^{-T}\GE\; (\Gt_1\,|\, \Gt_2\,|\, \Gn)^{-1}\qquad\hbox{weakly in}\quad (L^2(\Omega))^9,\cr}\right.\leqno(9.4)$$ where 
$$\GE=\pmatrix{
\displaystyle  S_3{\partial \GR\over \partial s_1} \Gn \cdot\GR\Gt_1+{\cal Z}_{11} & \displaystyle  S_3 {\partial \GR\over \partial s_1} \Gn \cdot\GR\Gt_2 +{\cal Z}_ {12}
 &\displaystyle  {1\over 2}\GR^T{\partial\overline{V} \over \partial S_3} \cdot \Gt_1\cr
 * & \displaystyle  S_3{\partial \GR\over \partial s_2} \Gn \cdot\GR\Gt_2+{\cal Z}_{22}   &\displaystyle  {1\over 2}\GR^T{\partial\overline{V} \over \partial S_3} \cdot \Gt_2 \cr
 * & *&  \displaystyle  \GR^T{\partial\overline{V} \over \partial S_3} \cdot \Gn \cr}$$ 
Now, recall that 
$${J_{2,\delta}(v_\delta)\over \delta^3}
=\int_{\Omega}{1\over \delta^2} W\Big({1\over 2}\Pi _\delta \big((\nabla_xv_\delta)^T\nabla_x v_\delta-\GI_3\big)\Big) \Pi_\delta\big(\det(\nabla\Phi)\big)
-{1\over \delta^3}\int_{{\cal Q}_\delta} f_{\kappa,\delta} \cdot  (v_\delta-I_d).\leqno(9.5)$$ 

In order to pass to the lim-inf in (9.5) we first notice that $\ds \det(\nabla\Phi)=\det (\Gt_1 | \Gt_2Ê| \Gn) +s_3\det\Big({\partial \Gn \over \partial s_1} |\Gt_2\; |\; \Gn\Big)+s_3\det\Big(\Gt_1 | {\partial \Gn \over \partial s_2} \; |\; \Gn\Big)+s_3^2\det\Big( {\partial \Gn \over \partial s_1} | {\partial \Gn \over \partial s_2} \; |\; \Gn\Big)$ so that indeed $\Pi_\delta\big(\det(\nabla\Phi)\big)$ strongly converges to $\det (\Gt_1 | \Gt_2Ê| \Gn) $ in $L^\infty(\Omega)$ as $\delta$ tends to $0$.

We  now consider the first term of the right hand side. Let $\varepsilon>0$ be fixed. Due to (6.2), there exists $\theta>0$ such that 
$$\forall E\in\GS_3,  \enskip   |||E|||\le \theta, \enskipÊ W(E)\ge Q(E)-\varepsilon |||E|||^2.\leqno (9.6)$$
We now use a similar argument given in [5]. Let us denote by $\chi_\delta^\theta$ the characteristic function of the set $A_\delta^\theta=\{s\in\Omega; |||\Pi _\delta \big((\nabla_xv_\delta)^T\nabla_x v_\delta-\GI_3\big)(s)||| \geq \theta\}$. Due to (9.2), we have
$$ {\rm meas}(A_\delta^\theta)\le C{\delta^2\over \theta^2} .\leqno(9.7)$$
Using the positive character of $W$, (9.2) and (9.6) give
$$\eqalign {
\int_{\Omega}{1\over \delta^2} \widehat{W}\big(\Pi _\delta (\nabla_xv_\delta)\big)|\Pi_\delta\big(\det(\nabla\Phi)\big)| &\ge \int_{\Omega}{1\over \delta^2} W\Big({1\over 2}\Pi _\delta \big((\nabla_xv_\delta)^T\nabla_x v_\delta-\GI_3\big)\Big)(1-\chi_\delta^\theta)\Pi_\delta\big(\det(\nabla\Phi)\big)\cr
&\ge  \int_{\Omega} Q\Big({1\over 2\delta}\Pi _\delta \big((\nabla_xv_\delta)^T\nabla_x v_\delta-\GI_3\big)(1-\chi_\delta^\theta)\Big)\Pi_\delta\big(\det(\nabla\Phi)\big)-C\varepsilon\cr}$$
In view of (9.7), the function $\chi_\delta^\theta$ converges  a.e. to $0$ as $\delta$ tends to $0$ while the weak limit of $\ds{1\over 2\delta}\Pi _\delta \big((\nabla_xv_\delta)^T\nabla_x v_\delta-\GI_3\big)(1-\chi_\delta^\theta)$ is given by (9.4). As a consequence  and also using the convergence of $\Pi_\delta\big(\det(\nabla\Phi)\big)$ obtained above, we have 
$$\liminf_{\delta\to 0} \int_{\Omega}{1\over \delta^2}\widehat{W}\big(\Pi _\delta (\nabla_xv_\delta)\big)\Pi_\delta\big(\det(\nabla\Phi)\big)\ge 
\int_{\Omega} Q\Big((\Gt_1\,|\, \Gt_2\,|\, \Gn)^{-T}\GE\; (\Gt_1\,|\, \Gt_2\,|\, \Gn)^{-1}\Big)\det (\Gt_1 | \Gt_2Ê| \Gn)-C\varepsilon.$$
As $\varepsilon$ is arbitrary, this gives
$$\liminf_{\delta\to 0} \int_{\Omega}{1\over \delta^2} \widehat{W}\big(\Pi _\delta (\nabla_xv_\delta)\big)\Pi_\delta\big(\det(\nabla\Phi)\big)
\ge \int_{\Omega} Q\Big((\Gt_1\,|\, \Gt_2\,|\, \Gn)^{-T}\GE\; (\Gt_1\,|\, \Gt_2\,|\, \Gn)^{-1}\Big)\det (\Gt_1 | \Gt_2Ê| \Gn).\leqno (9.8)$$
Using the convergences (9.4), it follows that
$$\lim_{\delta \to 0}\Big({1\over \delta^3}\int_{{\cal Q}_\delta} f_{2,\delta} \cdot  (v_\delta-I_d)\Big)={\cal L}({\cal V},\GR)$$
where ${\cal L}(\cdot,\cdot)$ is defined by (8.5). From (9.5), (9.8) and the above limit, we conclude that 
$$\liminf_{\delta\to 0}{m_{2,\delta}\over \delta^{3}}=\lim_{\delta\to 0}{J_{2,\delta}(v_\delta)\over \delta^3}\ge \int_{\Omega} Q\Big((\Gt_1\,|\, \Gt_2\,|\, \Gn)^{-T}\GE\; (\Gt_1\,|\, \Gt_2\,|\, \Gn)^{-1}\Big)\det (\Gt_1 | \Gt_2Ê| \Gn)-{\cal L}({\cal V},\GR).\leqno (9.9)$$ Proceeding as in the proof  of (8.7) in Section 8, we get
$$ \int_{\Omega} Q\Big((\Gt_1\,|\, \Gt_2\,|\, \Gn)^{-T}\GE\; (\Gt_1\,|\, \Gt_2\,|\, \Gn)^{-1}\Big)\det (\Gt_1 | \Gt_2Ê| \Gn)-{\cal L}({\cal V},\GR)\ge 
\min_{\Gv\in \D}{\cal J}_2(\Gv)=m^s_2.$$ Finally we have proved that $\ds m^s_2\le \liminf_{\delta\to 0}{m_{2,\delta}\over \delta^{3}}$.
\smallskip
\noindent{\ggras Step 2.} In this step we prove that $\ds m^s_2\ge  \limsup_{\delta\to 0}{m_{2,\delta}\over \delta^{3}}$. 
\smallskip
\noindent Let us now consider a minimizer $\Gv_0=\big({\cal V}_0,\GR_0,\overline{V}_0\big)\in \D$ of ${\cal J}_2$ and the  sequence $\Big(\big({\cal V}_\delta, \GR_\delta,\overline{V}_\delta\big)\Big)_{\delta>0}$ of  approximation of  $\Gv_0$ given by Lemma C constructed in the Appendix. The deformation $v_\delta$ is now defined by
$$v_{\delta}(s)={\cal V}_\delta(s_1,s_2)+s_3 \GR_\delta(s_1,s_2)\Gn(s_1,s_2)+\delta^2\overline{V}_\delta\Big(s_1,s_2,{s_3\over \delta}\Big),\quad \hbox{for}\  s\in \Omega_\delta.\leqno(9.10)$$ 
\noindent {\bf Step 2.1.} Estimate on $||\Pi_\delta\big(\nabla_xv_\delta-\GR_\delta\big)||_{(L^\infty(\Omega))^{3\times 3}}$ and $ ||dist \big(\nabla_x v_\delta,SO(3)\big)||_{L^\infty(\omega)} $.
\medskip
From (9.10) and trough simple calculations, we first have
$$\left\{\eqalign{
(\nabla_x v_\delta-\GR_\delta)\Gt_\alpha&={\partial{\cal V}_\delta\over \partial s_\alpha}-\GR_\delta\Gt_\alpha+s_3{\partial \GR_\delta\over \partial s_\alpha}\Gn+\delta^2 {\partial\overline{V}_\delta\over \partial s_\alpha}-(\nabla_x v_\delta-\GR_\delta)s_3{\partial\Gn\over \partial s_\alpha}\cr
(\nabla_x v_\delta-\GR_\delta)\Gn&=\delta {\partial\overline{V}_\delta\over \partial S_3},\cr}\right.\leqno(9.11)$$
then 
$$\eqalign{
&\Pi_\delta(\nabla_x v_\delta-\GR_\delta)\cdot \Pi_\delta(\nabla_s \Phi)\cr
=&\Big({\partial{\cal V}_\delta\over \partial s_1}-\GR_\delta\Gt_1+S_3\delta{\partial \GR_\delta\over \partial s_\alpha}\Gn+\delta^2{\partial \overline{V}_\delta\over \partial s_\alpha}\;|\; {\partial{\cal V}_\delta\over \partial s_2}-\GR_\delta\Gt_2+S_3\delta{\partial \GR_\delta\over \partial s_\alpha}\Gn+\delta^2{\partial \overline{V}_\delta\over \partial s_\alpha}\;|\; \delta {\partial\overline{V}_\delta\over \partial S_3}\Big).\cr}\leqno(9.12)$$
Thanks to (2.3) and the estimates of Lemma C in Appendix we obtain
$$||\Pi_\delta\big(\nabla_xv_\delta-\GR_\delta\big)||_{(L^\infty(\Omega))^{3\times 3}}\le {1\over 4}\leqno(9.13)$$ and  we  deduce that there exists a positive constant $C_0$ such that
$$||\Pi_\delta\big((\nabla_xv_\delta)^T\nabla_xv_\delta-\GI_3\big)||_{(L^\infty(\Omega))^{3\times 3}}\le C_0.\leqno(9.14)$$
Again using the estimates in Lemma C  we get
$$ ||dist \big(\nabla_x v_\delta,SO(3)\big)||_{L^\infty(\omega)} \le {1\over 2}$$ and then we obtain
$$\hbox{for a.e. }\enskip s\in \Omega_\delta \qquad \det\big(\nabla_x v_\delta(s)\big)>0.\leqno(9.15)$$

\noindent{\bf Step 2.2.} Strong limit of $\ds{1\over 2\delta}\Pi_\delta\Big((\nabla_xv_\delta)^T\nabla_x v_\delta-\GI_3\Big)$.
\medskip
Thanks to  the estimates and convergences of Lemma C  and (9.12)  we have
$$||\Pi_\delta\big(\nabla_xv_\delta-\GR_\delta\big)||_{(L^2(\Omega))^{3\times 3}}\le C\delta.\leqno(9.16)$$
We write the identity $(\nabla_xv_\delta)^T\nabla_x v_\delta-\GI_3=(\nabla_xv_\delta-\GR_\delta)^T\GR_\delta+\GR^T_\delta(\nabla_xv_\delta-\GR_\delta)+(\nabla_xv_\delta-\GR_\delta)^T(\nabla_xv_\delta-\GR_\delta)+(\GR_\delta-\GR)^T\GR_\delta+\GR^T(\GR_\delta-\GR)$. So, from (9.13) and (9.16) we get
$$ ||\Pi_\delta\big((\nabla_xv_\delta)^T\nabla_x v_\delta-\GI_3\big)||_{(L^2(\Omega))^{3\times 3}}\le C\delta.\leqno(9.17)$$
In view of (9.11), the strong convergences of Lemma C and  (9.16)  we deduce that
$$\left\{\eqalign{
{1\over \delta}\Pi_\delta\big((\nabla_x v_\delta-\GR)\Gt_\alpha\big)&\longrightarrow S_3{\partial \GR\over \partial s_\alpha}\Gn\quad \hbox{strongly in }\quad (L^2(\Omega))^3\cr
{1\over \delta}\Pi_\delta\big((\nabla_x v_\delta-\GR)\Gn\big)&\longrightarrow  {\partial\overline{V}\over \partial S_3}\cdot\Gn\quad \hbox{strongly in }\quad (L^2(\Omega))^3\cr}\right.\leqno(9.18)$$ Now thanks (9.13)  and the strong convergences (9.18) we obtain
$${1\over \sqrt\delta}\Pi_\delta (\nabla_x v_\delta-\GR)\longrightarrow 0\quad \hbox{strongly in }\quad (L^4(\Omega))^3$$ and then using again Lemma C, (9.18) and the above decomposition of $(\nabla_x v_\delta)^T\nabla_x v_\delta-\GI_3$, we get
$${1\over 2\delta}\Pi_\delta\big((\nabla_x v_\delta)^T\nabla_x v_\delta-\GI_3\big)\longrightarrow (\Gt_1 |\Gt_2 |\Gn)^{-T}\GE(\Gv_0)(\Gt_1 |\Gt_2 |\Gn)^{-1}\quad \hbox{strongly in }\quad (L^2(\Omega))^{3\times 3},\leqno(9.19)$$ where $\GE(\Gv_0)$ is given by (8.2).
\medskip
\noindent{\bf Step 2.3.} Let $\varepsilon $ be a fixed positive constant and  let   $\theta$ given by (7.2). We denote $\chi^\theta_\delta$ the characteristic function of the set $A^\theta_\delta=\{s\in\Omega; |||\Pi _\delta \big((\nabla_xv_\delta)^T\nabla_x v_\delta-\GI_3\big)(s)||| \ge \theta\}$. Due to (9.17), we have
$$ {\rm meas}(A^\theta_\delta)\le C{\delta^2\over \theta^2}\leqno(9.20)$$ and from (9.15) we have  $\det\big( \nabla_x v_\delta(s)\big)>0$ for a. e. $s\in \Omega_\delta$. Due to (6.2), (6.4) and (9.19)   we deduce that
$$\eqalign{
\limsup_{\delta\to 0} \int_{\Omega}{1\over \delta^2}(1-\chi^\theta_\delta) \, \widehat{W}\big(\Pi _\delta \big(\nabla_xv_\delta\big)\big)\Pi_\delta\big(\det(\nabla\Phi)\big)  &\le
\int_{\Omega} Q\Big((\Gt_1\,|\, \Gt_2\,|\, \Gn)^{-T}\GE(\Gv_0) (\Gt_1\,|\, \Gt_2\,|\, \Gn)^{-1}\Big)\det (\Gt_1 | \Gt_2Ê| \Gn)\cr
&+\varepsilon\int_{\Omega} |||(\Gt_1\,|\, \Gt_2\,|\, \Gn)^{-T}\GE(\Gv_0)  (\Gt_1\,|\, \Gt_2\,|\, \Gn)^{-1}|||^2\det (\Gt_1 | \Gt_2Ê| \Gn)\cr}$$  where $\GE(\Gv_0) $ is given by (8.2). Notice that there exists  a positive constant $C_1$ such that for all $ E\in \GS_3$ satisfying $\theta\le |||E|||\le C_0 $ we have
$$W(E)\le C_1|||E|||.$$
Thanks to  (6.3), (6.4), (9.17), the strong convergence (9.19)  and the weak convergence $\ds{1\over \delta}\chi^\theta_\delta\rightharpoonup 0 $ in $L^2(\Omega)$  we obtain
$$\lim_{\delta\to 0} \int_{\Omega}{1\over \delta^2}\chi^\theta_\delta \, \widehat{W}\big(\Pi _\delta (\nabla_xv_\delta)\big)\Pi_\delta\big(\det(\nabla\Phi) \big)\le C_1 \lim_{\delta\to 0} \int_{\Omega}{1\over \delta}\chi^\theta_\delta \, |||{1\over 2\delta}\Pi_\delta\big((\nabla_x v_\delta)^T\nabla_x v_\delta-\GI_3\big)|||\, \Pi_\delta\big(\det(\nabla\Phi)\big)=0$$
Hence for any $\varepsilon>0$ we get
$$\eqalign{
\limsup_{\delta\to 0} \int_{\Omega}{1\over \delta^2} \, \widehat{W}\big(\Pi _\delta (\nabla_xv_\delta)\big)\Pi_\delta\big(\det(\nabla\Phi)\big) &\le
\int_{\Omega} Q\Big((\Gt_1\,|\, \Gt_2\,|\, \Gn)^{-T}\GE (\Gv_0)(\Gt_1\,|\, \Gt_2\,|\, \Gn)^{-1}\Big)\det (\Gt_1 | \Gt_2Ê| \Gn)\cr
&+\varepsilon\int_{\Omega} |||(\Gt_1\,|\, \Gt_2\,|\, \Gn)^{-T}\GE (\Gv_0)(\Gt_1\,|\, \Gt_2\,|\, \Gn)^{-1}|||^2\det (\Gt_1 | \Gt_2Ê| \Gn)\cr}$$
Finally 
$$\limsup_{\delta\to 0} \int_{\Omega}{1\over \delta^2} \, \widehat{W}\big(\Pi _\delta (\nabla_xv_\delta)\big)\Pi_\delta\big(\det(\nabla\Phi)\big) \le
\int_{\Omega} Q\Big((\Gt_1\,|\, \Gt_2\,|\, \Gn)^{-T}\GE(\Gv_0) (\Gt_1\,|\, \Gt_2\,|\, \Gn)^{-1}\Big)\det (\Gt_1 | \Gt_2Ê| \Gn).\leqno(9.21)$$
As far as the contribution of the applied forces is concerned, we use the convergences of Lemma C  to obtain
$$\lim_{\delta \to 0}\Big({1\over \delta^3}\int_{{\cal Q}_\delta} f_{2,\delta} \cdot (v_\delta-I_d)\Big)={\cal L}({\cal V},\GR).\leqno(9.22)$$
From (9.21) and (9.22), we conclude that 
$$\eqalign{
\limsup_{\delta\to 0}{J_{2,\delta}(v_\delta)\over \delta^3} &\le \int_{\Omega} Q\Big((\Gt_1\,|\, \Gt_2\,|\, \Gn)^{-T}\GE(\Gv_0) (\Gt_1\,|\, \Gt_2\,|\, \Gn)^{-1}\Big)\det (\Gt_1 | \Gt_2Ê| \Gn)-{\cal L}({\cal V},\GR)={\cal J}_{2}(\Gv_0)=m^s_2.\cr}$$ Then we get $\ds \limsup_{\delta\to 0}{m_{2,\delta}\over \delta^{3}}\le m^s_2$.\fin
\noindent {\Ggras 10. Alternative formulations of the minima $m^s_{\kappa,\delta}$ and $m^s_2$.} 
\smallskip

In the following theorem we characterize the minimum of the functional $J^s_{\kappa,\delta}(\cdot) $ over $\D_{\delta,\gamma_0}$, respectively ${\cal J}_2$ over $\D$, as the minima of two functionals which depend on the mid-surface deformation ${\cal V}$ and on the matrix $\GR$ which gives the rotation of the fibers.

The first  theorem of this section shows that the variable  $\overline{v}$  can be eliminated in the minimization problem (7.9). 
\smallskip
\noindent We set
$$\eqalign{
\E=\Bigl\{ &\big({\cal V},\GR\big)\in (H^1(\omega))^3\times (H^1(\omega))^{3\times 3}\;\; |\;\; 
{\cal V}=\phi,\quad \GR=\GI_3\enskip\hbox{on}\; \gamma_0,\cr
&\quad\GR(s_1,s_2)\in SO(3)\quad  \hbox{for a.e.}\enskip (s_1,s_2)\in \omega\Big\}.\cr}$$ We recall (see (5.3)) that for all $({\cal V},\GR)\in \E$ we have set
$$\eqalign{
{\cal Z}_{\alpha\beta}&={1\over 2} \Big[\Big({\partial{\cal V}\over \partial s_\alpha }-\GR\Gt_\alpha\Big)\cdot\GR\Gt_\beta+\Big({\partial{\cal V}\over \partial s_\beta}-\GR\Gt_\beta\Big)\cdot\GR\Gt_\alpha\Big],\quad  {\cal Z}_{3\alpha}= {\partial{\cal V}\over \partial s_\alpha }\cdot\GR\Gn\cr
\Gamma_{\alpha\beta}(\GR)&={1\over 2}\Big\{{\partial \GR\over \partial s_\alpha} \Gn \cdot \GR\Gt_\beta+{\partial \GR\over \partial s_\beta} \Gn \cdot \GR\Gt_\alpha\Big\}.\cr}$$
 \noindent{\ggras  Theorem  10.1. }{\it Let $\Gv_\delta=({\cal V}_\delta, \GR_\delta , \overline{V}_\delta)\in \D_{\delta,\gamma_0}$ such that $\ds m^s_{\kappa,\delta}=J^s_{\kappa,\delta}(\Gv_\delta) =\min_{v\in \D_{\delta,\gamma_0}}J^s_{\kappa,\delta}(\Gv)$. We have
 $$m^s_{\kappa,\delta}={\cal F}^s_{\kappa,\delta}\big({\cal V}_\delta,\GR_\delta\big)=\min_{\ds({\cal V},\GR)\in\E}{\cal F}^s_{\kappa,\delta}\big({\cal V},\GR\big) \leqno(10.1) $$ where 
$$\left\{\eqalign{
{\cal F}^s_{\kappa,\delta}\big({\cal V},\GR\big)&=\delta^3\int_\omega a_{\alpha\beta\alpha^{'}\beta^{'}}\Gamma_{\alpha\beta}(\GR)\Gamma_{\alpha^{'}\beta^{'}}(\GR)+  \delta\int_\omega b_{i\alpha i^{'}\alpha^{'}}{\cal Z}_{i\alpha}{\cal Z}_{i^{'}\alpha^{'}}\cr
&+\delta^3\Big\|{\partial \GR\over \partial s_1} \Gt_2-{\partial \GR\over \partial s_2} \Gt_1\Big\|^2_{(L^2(\omega))^{3}}+\delta\Big\|{\partial{\cal V}\over \partial s_1 }\cdot\GR\Gt_2- {\partial{\cal V}\over \partial s_2}\cdot\GR\Gt_1\Big\|^2_{L^2(\omega)}-\delta^{\kappa^{'}+1}{\cal L}({\cal V}, \GR).\cr}\right.\leqno(10.2)$$ The $a_{\alpha\beta\alpha^{'}\beta^{'}}$ and $ b_{i\alpha i^{'}\alpha^{'}}$ are constants which depend only of the quadratic form $Q$ and the vectors $(\Gt_1,\Gt_2,\Gn)$.}
\smallskip
\noindent{\ggras Proof. } We have
$$\eqalign{
 m^s_{\kappa,\delta}=\min_{\Gv\in \D_{1,\gamma_0}} \Big[& \delta\int_\Omega Q\Big((\Gt_1\,|\, \Gt_2\,|\, \Gn)^{-T}\Pi_\delta( \widehat{E}(\Gv))(\Gt_1\,|\, \Gt_2\,|\, \Gn)^{-1}\Big)\det (\Gt_1 | \Gt_2Ê| \Gn)ds_1ds_2dS_3\cr
&+\delta^3\Big\|{\partial \GR\over \partial s_1} \Gt_2-{\partial \GR\over \partial s_2} \Gt_1\Big\|^2_{(L^2(\omega))^{3}}+\delta\Big\|{\partial{\cal V}\over \partial s_1 }\cdot\GR\Gt_2- {\partial{\cal V}\over \partial s_2}\cdot\GR\Gt_1\Big\|^2_{L^2(\omega)}-\delta^{\kappa^{'}+1}{\cal L}({\cal V}, \GR)\Big].\cr}$$  
In order to eliminate $ \overline{v}$, we first fix  $\big({\cal V} ,\GR \big)\in \E$.  We set 
$$\int_{-1}^1 Q\Big((\Gt_1\,|\, \Gt_2\,|\, \Gn)^{-T}\Pi_\delta( \widehat{E}(\Gv))(\Gt_1\,|\, \Gt_2\,|\, \Gn)^{-1}\Big)\det (\Gt_1 | \Gt_2Ê| \Gn)dS_3={\cal Q}_m(\Ga,\Gb,\Gc+\Gd)$$ where 
$$\Ga=\delta\pmatrix{\Gamma_{11}(\GR)\cr \Gamma_{12}(\GR)\cr \Gamma_{22}(\GR)\cr},\quad \Gb=\pmatrix{{\cal Z}_{11}\cr {\cal Z}_{12}\cr {\cal Z}_{22}\cr},\quad \Gd=\pmatrix{{\cal Z}_{31}\cr{\cal Z}_{32}\cr 0\cr},\qquad  \Gc={1\over \delta} \pmatrix{\ds {1\over 2}\GR^T {\partial\Pi_\delta(\overline{v}) \over \partial S_3} \cdot \Gt_1\cr \ds {1\over 2}\GR^T{\partial\Pi_\delta(\overline{v})\over \partial S_3} \cdot \Gt_2\cr\ds  \GR^T{\partial\Pi_\delta(\overline{v})\over \partial S_3} \cdot \Gn\cr}$$  and we apply Lemma B in Appendix  to obtain  the theorem. \fin
\smallskip
The next theorem is similar to Theorem 10.1 for the limit energy and the minimization problem (8.3). We set
$$\eqalign{
\E_{lim}=\Bigl\{ &\big({\cal V},\GR\big)\in \E\;\; |\;\; 
 {\partial{\cal V} \over \partial s_\alpha}=\GR \Gt_\alpha,\;\; \alpha=1,2\Big\}.\cr}$$
\noindent{\ggras  Theorem  10.2. }{\it Let $\Gv_0=({\cal V}_0, \GR_0 , \overline{V}_0)\in \D$ such that $\ds m^s_2={\cal J}_2(\Gv_0) =\min_{v\in \D}{\cal J}_2(\Gv)$. We have
$$\ds m^s_2={\cal F}_{2}\big({\cal V}_0,\GR_0\big)=\min_{ \ds ({\cal V},\GR)\in\E_{lim}}{\cal F}_{2}
\big({\cal V},\GR\big) \leqno(10.3) $$ where  
$${\cal F}_{2}\big({\cal V},\GR\big)=\int_\omega a_{\alpha\beta\alpha^{'}\beta^{'}}\Big({\partial \GR\over \partial s_\alpha} \Gn \cdot\GR\Gt_\beta\Big)\Big({\partial \GR\over \partial s_\alpha^{'}} \Gn \cdot\GR\Gt_{\beta^{'}}\Big)-{\cal L}({\cal V},\GR)\leqno(10.4)$$ The $a_{\alpha\beta\alpha^{'}\beta^{'}}$ are the same constants as the one in Theorem 10.1.}
\smallskip
\noindent{\ggras Proof. } We proceed as in Theorem 10.1. In order to eliminate $ \overline{V}$, we fix  $\big({\cal V} ,\GR \big)\in \E_{lim}$  and we minimize the functional ${\cal J}_{2}\big({\cal V} ,\GR ,\cdot\big)$ over the space $\GV$.  Thanks to Lemma B in Appendix we obtain the minimum with respect to $\overline{V}$ and then the new characterization of the minimum $m^s_2$.\fin

Of course, for all $({\cal V},\GR)\in \E_{lim}$, we get
$${\cal F}_{2,\delta}^s\big({\cal V},\GR\big)=\delta^3{\cal F}_{2}\big({\cal V},\GR\big).$$

  Let us give the explicit expression of the limit energies ${\cal F}^s_{\kappa,\delta}$ and ${\cal F}_{2}$ in the case where $S$ is a developable surface such that the parametrization $\phi$ is locally isometric
  $$\forall (s_1,s_2) \in \overline{\omega} \qquad ||\Gt_\alpha(s_1,s_2)||_2 =1\qquad \Gt_1(s_1,s_2)\cdot \Gt_2(s_1,s_2)=0.$$
  We consider a  St Venant-Kirchhoff's law for which we have
$$\widehat{W}(F)=\left\{\eqalign{
&{\lambda\over 8}\big(tr(F^TF-\GI_3) \big)^2+{\mu\over 4}tr\big((F^TF-\GI_3)^2\big)\quad\hbox{if}\quad \det(F)>0\cr
&+\infty\qquad \hbox{if}\qquad \det(F)\le 0,\cr}\right.$$  so that $Q=W=W^s$.
\smallskip
\noindent{\ggras  Expression of  ${\cal F}^s_{\kappa,\delta}$.} For any $\Gv=\big({\cal V},\GR,\overline{v}\big)\in \D_{\delta,\gamma_0}$, the expression (7.5) gives
$$\left\{\eqalign{
J^s_{\kappa,\delta}(\Gv)=\delta& \int_{\Omega} \Big[{\lambda\over 2}\big(tr(\widehat{E}(\Gv)) \big)^2+\mu\, tr\big((\widehat{E}(\Gv))^2\big)\Big]+ \delta^3\Big\|{\partial \GR\over \partial s_1} \Gt_2-{\partial \GR\over \partial s_2} \Gt_1\Big\|^2_{(L^2(\omega))^{3}}\cr
&+\delta \Big\|{\partial{\cal V}\over \partial s_1 }\cdot\GR\Gt_2- {\partial{\cal V}\over \partial s_2}\cdot\GR\Gt_1\Big\|^2_{L^2(\omega)} - \delta^{\kappa^{'}+1}{\cal L}({\cal V}, \GR).\cr}\right.\leqno(10.5)$$ where $\widehat{E}(\Gv)$ is defined by (5.3). It follows that the elimination of $\overline{V}$ in Theorem 10.1 gives the partial derivatives of $\overline{V}$ with respect to $S_3$
$$\pmatrix{
\ds{\partial\overline{v}\over \partial s_3} (.,., s_3)\cdot\GR \Gt_1\cr\cr
\ds {\partial\overline{v}\over \partial s_3} (.,., s_3)\cdot\GR \Gt_2\cr\cr
\ds{\partial\overline{v}\over \partial s_3}(.,., s_3)\cdot\GR \Gn\cr}=\pmatrix{
\ds -{{\cal Z}_{31}\over \delta^2}\Big(\delta^2+{5\over 4}(s^2_3-\delta^2)\Big)\cr\cr
\ds -{{\cal Z}_{32}\over \delta^2}\Big(\delta^2+{5\over 4}(s^2_3-\delta^2)\Big)\cr\cr
\ds -{\nu\over 1-\nu}\Big(s_3\big[\Gamma_{11}(\GR)+\Gamma_{22}(\GR)\big]+\big[{\cal Z}_{11}+{\cal Z}_{22}\Big]\Big)\cr}\leqno(10.6)$$ and then 
$$\eqalign{
{\cal F}^s_{\kappa,\delta}\big({\cal V},\GR\big)=&{E\delta^3\over 3(1-\nu^2)}\int_\omega\Big[(1-\nu)\sum_{\alpha,\beta=1}^2\big(\Gamma_{\alpha\beta}(\GR)\big)^2+\nu\big(\Gamma_{11}(\GR)+\Gamma_{22}(\GR)\big)^2\Big]\cr
&+{E\delta\over (1-\nu^2)}\int_\omega\Big[(1-\nu)\sum_{\alpha,\beta=1}^2\big({\cal Z}_{\alpha\beta}\big)^2+\nu\big({\cal Z}_{11}+{\cal Z}_{22}\big)^2\Big]+{5E\delta\over 12(1+\nu)}\int_\omega\big({\cal Z}_{31}^2+{\cal Z}_{32}^2\big)\cr
&+ \delta^3\Big\|{\partial \GR\over \partial s_1} \Gt_2-{\partial \GR\over \partial s_2} \Gt_1\Big\|^2_{(L^2(\omega))^{3}}
+ \delta\Big\|{\partial{\cal V}\over \partial s_1 }\cdot\GR\Gt_2- {\partial{\cal V}\over \partial s_2}\cdot\GR\Gt_1\Big\|^2_{L^2(\omega)} - \delta^{\kappa^{'}+1}{\cal L}({\cal V}, \GR).\cr}$$
\smallskip
\noindent{\ggras  Expression of  ${\cal F}_{2}$.} For any $\Gv=\big({\cal V},\GR,\overline{V}\big)\in \D$, the expression (8.1) gives
$${\cal J}_{2}\big( \Gv\big)= \int_{\Omega} \Big[{\lambda\over 2
}\big(tr(\GE(\Gv)) \big)^2+\mu\, tr\big((\GE(\Gv))^2\big)\Big]-{\cal L}({\cal V},\GR)$$ where $\GE(\Gv)$ is defined by (8.2). It follows that the elimination of $\overline{V}$ in Theorem 10.2 is identical to that of standard linear elasticity (see [18]) hence we have 
$$\overline{V} (.,., S_3)=-{\nu\over 2(1-\nu)}\Big(S_3^2-{1\over 3}\Big)\Big[\Gamma_{11}(\GR)+\Gamma_{22}(\GR)\Big]\GR\Gn\leqno(10.7)$$ and then 
$${\cal F}_{2}\big({\cal V},\GR\big)={E\over 3(1-\nu^2)}\int_\omega\Big[(1-\nu)\sum_{\alpha,\beta=1}^2\big(\Gamma_{\alpha\beta}(\GR)\big)^2+\nu\big(\Gamma_{11}(\GR)+\Gamma_{22}(\GR)\big)^2\Big]-{\cal L}({\cal V},\GR).$$
{\ggras Remark 10.1.} In the case of a St-Venant-Kirchhoff material  a classical energy argument show that if $\big(v_\delta\big)_{0<\delta\le \delta_0}$ is a sequence such that 
$$m^s_2=\lim_{\delta\to 0}{J_{2,\delta}(v_\delta)\over \delta^{3}},$$
then there exists a subsequence and $({\cal V}_0,\GR_0)\in \E$, which is a solution of Problem (10.3), such that  the sequence of the Green-St Venant's deformation tensors satisfies
 $${1\over 2\delta}\Pi _\delta \big((\nabla_xv_\delta)^T\nabla_x v_\delta-\GI_3\big) \longrightarrow   (\Gt_1\,|\, \Gt_2\,|\, \Gn)^{-T}\GE(\Gv_0) (\Gt_1\,|\, \Gt_2\,|\, \Gn)^{-T}  \qquad\hbox{strongly in}\quad (L^2(\Omega))^{3\times 3},$$
 where $\GE(\Gv_0)$ is defined in (8.2) with $\overline{V}_0$ given by (10.7) (replacing $\GR$ by $\GR_0$). 
\smallskip
\noindent{\ggras Remark 10.2.}   It is well known that the constraint $\ds{\partial{\cal V} \over \partial s_1}=\GR \Gt_1$ and  $\ds{\partial{\cal V} \over \partial s_2}=\GR \Gt_2$ together the boundary conditions are  strong limitations
  on the possible deformation for the limit 2d shell. Actually for a plate or as soon as $S$ is a developable surface, the configuration after deformation must also be  a developable surface. In the general case, it is an open problem to know if the set $\E_{lim}$ contains other deformations than  identity mapping or very special isometries (as for example symetries).
\medskip

\noindent{\Ggras Appendix.}
\medskip
\noindent{\Ggras Lemma A. }{\it Let ${\cal Q}_m$ be the  positive definite quadratic form defined on the space $\R^3\times \R^3\times \big(L^2(-1,1)\big)^3$ by
$$\forall(\Ga,\Gb,\Gc)\in \R^3\times \R^3\times \big(L^2(-1,1)\big)^3,\qquad {\cal Q}_m(\Ga,\Gb,\Gc)=\int_{-1}^1\GA(S_3)\pmatrix{S_3\Ga_1+\Gb_1\cr S_3\Ga_2+\Gb_2\cr S_3\Ga_3+\Gb_3\cr \Gc_1(S_3)\cr \Gc_2(S_3)\cr \Gc_3(S_3)\cr}\cdot\pmatrix{S_3\Ga_1+\Gb_1\cr S_3\Ga_2+\Gb_2\cr S_3\Ga_3+\Gb_3\cr \Gc_1(S_3)\cr \Gc_2(S_3)\cr \Gc_3(S_3)\cr}dS_3$$ where $\GA(S_3)$ is a symmetric positive definite $6\times 6$ matrix satisfying
$$ \GA(S_3)=\GA(-S_3)\qquad \hbox{for a.e. }S_3\in ]-1,1[\leqno(A.1)$$ and moreover there exists a positive constant $c$ such that
$$ \forall \xi\in \R^6,\qquad \GA(S_3)\xi\cdot\xi\ge c|\xi|^2\quad \qquad \hbox{for a.e. }S_3\in ]-1,1[.\leqno(A.2)$$ For all $\Ga\in \R^3$, we have
$$\min_{(\Gb,\Gc)\in  \R^3\times  (L^2(-1,1))^3 } {\cal Q}_m(\Ga,\Gb,\Gc)=\min_{\Gc\in  \GL_2} {\cal Q}_m(\Ga,0,\Gc)$$ where
$$\GL_2=\Big\{\Gc\in  (L^2(-1,1))^3\;|\; \int_{-1}^1 \Gc_\alpha(S_3)(S_3^2-1)dS_3=0,\enskip \alpha\in\{1,2\}\Big\}.$$}
\noindent{\Ggras Proof. } We write
$$\GA(S_3)=\pmatrix{ \GA_1(S_3) & \vdots & \GA_2(S_3) \cr 
\cdots & & \cdots \cr \GA_2^T(S_3) & \vdots &  \GA_3(S_3)\cr}$$ where for a.e. $S_3\in ]-1,1[$, $\GA_1(S_3)$ and $\GA_3(S_3)$ are symmetric positive definite $3\times 3$ matrices. The both minimum are obtained with 
$$\Gc_0(S_3) =-\GA_3^{-1}(S_3) \GA_2^T(S_3)S_3\Ga,\qquad  \Gb_0=0.$$ We have
$${\cal Q}_m(\Ga,0,\Gc_0)=\Big(\int_{-1}^1 S^2_3\big(\GA_1(S_3)-\GA_2(S_3)\GA_3^{-1}(S_3) \GA_2^T(S_3)\big)dS_3\Big)\Ga\cdot \Ga.\leqno(A.3)$$\fin
In the following lemma we use the same notation as in Lemma A.
\smallskip
\noindent{\Ggras Lemma B. }{\it Let $\Ga$, $\Gb$  be two fixed vectors in $\R^3$ and let $\Gd$ be a fixed vector in $\R^2\times\{0\}$. We have
$$\min_{\Gc\in \GL_2} {\cal Q}_m(\Ga,\Gb,\Gc+\Gd)=\Big( \int_{-1}^1 S^2_3\big[\GA_1(S_3)-\GA_2(S_3)\GA_3^{-1}(S_3)\GA^T_2(S_3)\big]\Big)\Ga\cdot\Ga+Q^{'}_m\big(\Gb,\Gd\big)\leqno(B.1)$$  where $Q^{'}_m$ is a positive definite quadratic form  which depends only on the matrix $\GA$. }
\smallskip
\noindent{\ggras Proof. }  Through  solving a simple variational problem, we find that the minimum of the functional $\Gc\longmapsto Q_m(\Ga,\Gb,\Gc+\Gd)$ over the space $\GL_2$ is obtained with 
$$\Gc(S_3)=-\Gd-  \GA_3^{-1}(S_3)\GA^T_2(S_3)(S_3\Ga+\Gb)+(S^2_3-1)\GA_3^{-1}(S_3)\Ge$$ where  $\Ge\in \R^2\times\{0\}$
$$  \Ge=e_1 \Ge_1+e_2 \Ge_2,\qquad \Ge_1=\pmatrix{1 \cr 0\cr 0\cr},\;  \Ge_2=\pmatrix{0 \cr 1\cr 0\cr}$$ is the solution of  the system
$$\Big[{4\over 3}\Gd-\Big(\int_{-1}^1  (S_3^2-1)\GA_3^{-1}(S_3)\GA^T_2(S_3)dS_3\Big)\Gb+\Big(\int_{-1}^1(S^2_3-1)^2\GA_3^{-1}(S_3)dS_3\Big)\Ge\Big]\cdot\Ge_\alpha=0,\quad \alpha=1,2.$$ 
  Notice that the matrix $\ds \int_{-1}^1(S^2_3-1)^2\GA_3^{-1}(S_3)dS_3$ is a $3\times 3$ symmetric positive definite matrix. Replacing $\Gc$ and $\Ge$ by their  values   we obtain $(B.1)$.\fin
\noindent {\Ggras Lemma C.}
{\it Let  $\big({\cal V},\GR,\overline{V}\big)$ be in $\GD_2$, there exists  a sequence  $\Big(\big({\cal V}_\delta, \GR_\delta,\overline{V}_\delta\big)\Big)_{\delta>0}$ of   $(W^{2,\infty}(\omega))^3\times(W^{1,\infty}(\omega))^{3\times 3}\times(W^{1,\infty}(\Omega))^3$ such that
$${\cal V}_\delta=\phi, \quad   \GR_\delta= \GI_3 \quad  \hbox{ on }\; \gamma_0,\qquad
\overline{V}_\delta=0,\; \hbox{ on }\; \gamma_0\times ]-1,1[,\leqno(C.1)$$
with
$$\left\{\eqalign{
 {\cal V}_\delta&\longrightarrow {\cal V}\quad \hbox{strongly in }\quad (H^2(\omega))^3\cr
 \GR_\delta&\longrightarrow \GR\quad \hbox{strongly in }\quad (H^1(\omega))^{3\times 3}\cr
 {1\over \delta}(\GR_\delta&-\GR)\longrightarrow 0\quad \hbox{strongly in }\quad (L^2(\omega))^{3\times 3}\cr
  {1\over \delta}\Big({\partial {\cal V}_\delta \over \partial s_\alpha}&- \GR_\delta\Gt_\alpha\Big)\longrightarrow 0\quad \hbox{strongly in }\quad (L^2(\omega))^{3}\cr
\overline{V}_\delta&\longrightarrow \overline{V}\quad \hbox{strongly in }\quad (L^2(\omega;H^1((-1,1)))^{3},\cr
\delta{\partial \overline{V}_\delta\over \partial s_\alpha}&\longrightarrow 0\quad \hbox{strongly in }\quad (L^2(\Omega))^{3},\cr}\right.\leqno(C.2)$$ and moreover
 $$\left\{\eqalign{
  & ||dist \big(\GR_\delta,SO(3)\big)||_{L^\infty(\omega)} \le {1\over 8},\qquad  \Big\|{\partial{\cal V}_\delta\over \partial s_\alpha}-\GR_\delta\Gt_\alpha\Big\|_{(L^\infty(\omega))^3} \le{1\over 8},\cr
 &  || \GR_\delta||^2_{(W^{1,\infty}(\omega))^{3\times 3}}+ ||\overline{V}_\delta ||^2_{((W^{1,\infty}(\Omega))^3} \le{1\over (4c^{'}_1\delta)^2}. \cr}\right.\leqno (C.3)$$ The constant $c^{'}_1$ is given by (2.3).}
\medskip
\noindent {\Ggras Proof.}  For $h>0$ small enough, consider a ${\cal C}^\infty_0(\R^2)$-function $ \psi_h$ such that $0\le \psi_h\le1$
$$\left \{\eqalign{
& \psi_h(s_1,s_2)=1\  \hbox{if } dist\big((s_1,s_2),\gamma_0\big)\le h\cr
& \psi_h(s_1,s_2)=0 \ \hbox{if } dist\big((s_1,s_2),\gamma_0\big)\ge 2h.\cr}\right.$$
Indeed we can assume that 

$$||\psi_h||_{W^{1,\infty}(\R^2)} \le {C\over h}, \qquad ||\psi_h||_{W^{2,\infty}(\R^2)} \le {C\over h^2}.\leqno (C.4)$$
Since $\omega$ is bounded with a Lipschitz boundary, we first extend the fields  ${\cal V}$ and  $\GR_\Gn=\GR\Gn$  into two fields of  $(H^2(\R^2))^3$ and  $(H^1(\R^2))^3$ (and we use the same notations for these extentions). We  define the $3\times3$ matrix field $\GR^{'}\in (H^1(\R^2))^{3\times 3}$ by the formula
$$\GR^{'}=\Big({\partial {\cal V}\over \partial s_1} | {\partial {\cal V}\over \partial s_2} | \GR_{\Gn}\Big)  \Big(\Gt_1 | \Gt_2 | \Gn\Big)^{-1}. \leqno(C.5)$$
By construction we have  $\displaystyle{\partial {\cal V}\over \partial s_\alpha}=\GR^{'}\Gt_\alpha$ in $\R^2$ and $\GR^{'}=\GR$ in $\omega$. At least, we introduce below the  approximations ${\cal V}_h$ and $\GR_h$ of ${\cal V}$ and $\GR$ as  restrictions to $\overline{\omega}$ of the following fields defined into $\R^2$:
$$\left\{\eqalign{
{\cal V}^{'}_h(s_1,s_2)&={1\over 9\pi h^2}\int_{B(0,3h)} {\cal V}(s_1+t_1,s_2+t_2)dt_1dt_2,\cr
\GR^{'}_h(s_1,s_2)&={1\over 9\pi h^2}\int_{B(0,3h)}\GR^{'}(s_1+t_1,s_2+t_2)dt_1dt_2,\cr}\right.\quad \hbox{a.e. } (s_1,s_2)\in \R^2\leqno(C.6)$$ and
$${\cal V}_h=\phi \psi_h +{\cal V}^{'}_h(1-\psi_h),\qquad \GR_h=\GI_3\psi_h +\GR^{'}_h (1-\psi_h ),\quad \hbox{in}\enskip \omega.\leqno(C.7)$$
Notice that we have
$$\eqalign{
&{\cal V}^{'}_h\in (W^{2,\infty}(\R^2))^3,\qquad \GR^{'}_h\in (W^{1,\infty}(\R^2))^{3\times3},\cr
&{\cal V}_h\in (W^{2,\infty}(\omega))^3,\qquad \GR_h\in (W^{1,\infty}(\omega))^{3\times3},\qquad {\cal V}_h=\phi,\quad \GR_h=\GI_3\enskip\hbox{on}\enskip\gamma_0.\cr}\leqno(C.8)$$
 Due to the definition (C.5) of  $\GR^{'}$  and  in view of (C.6)  we have
 $$\left\{\eqalign{
{\cal V}^{'}_h&\longrightarrow {\cal V} \quad \hbox{strongly in } (H^2(\R^2))^3,\cr
 \GR^{'}_h&\longrightarrow \GR^{'} \quad \hbox{strongly in } (H^1(\R^2))^{3\times 3}\cr
}\right.\leqno(C.9) $$ and thus using estimates (C.4)
$$\left\{\eqalign{
{\cal V}_h&\longrightarrow {\cal V} \quad \hbox{strongly in } (H^2(\omega))^3,\cr
 \GR_h&\longrightarrow \GR \quad \hbox{strongly in } (H^1(\omega))^{3\times 3}\cr
}\right.\leqno(C.10) $$ 
Moreover using again (C.6) and the fact that $\GR^{'}-\GR_h$ strongly converges to $0$ in $(H^1(\R^2))^{3\times 3}$ we deduce that 
$${1\over h}\big(\GR^{'}_h-\GR^{'})\longrightarrow 0 \quad \hbox{strongly in } (L^2(\R^2))^{3\times 3}$$
and then together with (C.4), (C.5),  (C.7) and (C.10) we get
$$\eqalign{
{1\over h}\big(\GR_h-\GR\big)&\longrightarrow 0 \quad \hbox{strongly in } (L^2(\omega))^{3\times 3},\cr
{1\over h}\big({\partial{\cal V}_h\over \partial s_\alpha}-\GR_h\Gt_\alpha\big)&\longrightarrow 0 \quad \hbox{strongly in } (L^2(\omega))^{3}.\cr}$$ 
We now turn to the estimate of the distance between $\GR_h(s_1,s_2)$ and $SO(3)$ for a.e. $(s_1,s_2)\in\omega$. We apply  the PoincarŽ-Wirtinger's inequality to the function $(u_1,u_2)\longmapsto \GR^{'}(u_1,u_2)$ in the ball $B((s_1,s_2),3h)$. We obtain
$$\int_{B((s_1,s_2),3h)} |||\GR^{'}(u_1,u_2)-\GR^{'}_h(s_1,s_2)|||^2du_1du_2\le Ch^2||\nabla \GR^{'}||^2_{(L^2(B((s_1,s_2), 3h)))^3}$$ where $C$ is the PoincarŽ-Wirtinger's constant for a ball. Since the open set $\omega$ is boundy with a Lipschitz boundary, there exists a positive constant $c(\omega)$, which depends only on $\omega$, such that
$$|\big(B((s_1,s_2), 3h)\setminus B((s_1,s_2), 2h)\big)\cap \omega|\ge c (\omega)h^2.$$ Setting $m_h(s_1,s_2)$ the essential infimum of the function $(u_1,u_2)\mapsto |||\GR(u_1,u_2)-\GR^{'}_h(s_1,s_2)|||$ into the set $\big(B((s_1,s_2), 3h)\setminus B((s_1,s_2), 2h)\big)\cap \omega$, we then obtain
$$c(\omega) h^2m_h(s_1,s_2)^2\le Ch^2||\nabla \GR^{'}||^2_{(L^2(B((s_1,s_2), 3h)))^3}$$ Hence, thanks to the strong convergence of $\GR^{'}_h$ given by (C.9), the above inequality shows that there exists $h^{'}_0$ which does not depend on $(s_1,s_2)\in \overline{\omega}$ such that for any $h\le h^{'}_0$
$$dist(\GR^{'}_h(s_1,s_2) , SO(3))\le {1/8}\qquad \hbox{for any } (s_1,s_2)\in \overline{\omega}.$$ 
Now, 

$\bullet$ in the case $dist\big((s_1,s_2),\gamma_0\big)>2h$, $(s_1,s_2)\in \omega$,  by definition of $\GR_h$ and thanks to the above inequality we have $dist(\GR_h(s_1,s_2) , SO(3))\le {1/8}$,

$\bullet$ in the case  $dist\big((s_1,s_2),\gamma_0\big)< h$, $(s_1,s_2)\in \omega$,  by definition of $\GR_h$ we have $\GR_h(s_1,s_2)=\GI_3$ and then $dist(\GR_h(s_1,s_2) , SO(3))=0$,

$\bullet$ in the case  $h\le dist\big((s_1,s_2),\gamma_0\big)\le 2h$, $(s_1,s_2)\in \omega$,  due to the fact that $\GR^{'}=\GI_3$ onto $\gamma_0$, firstly we have
$$||\GR^{'}-\GI_3||^2_{(L^2(\omega_{6h,\gamma_0}))^{3\times 3}}\le Ch^2||\nabla \GR^{'}||^2_{(L^2(\omega_{6h,\gamma_0}))^{3\times 3}}$$ where $\omega_{kh,\gamma_0}=\{(s_1,s_2)\in \R^2\;|\; dist\big((s_1,s_2),\gamma_0\big)\le kh\}$, $k\in \N^*$. Hence
$$||\GR^{'}_h-\GI_3||^2_{(L^2(\omega_{3h,\gamma_0}))^{3\times 3}}\le Ch^2||\nabla \GR^{'}||^2_{(L^2(\omega_{6h,\gamma_0}))^{3\times 3}}.$$ The constants depend only on $\partial \omega$.
\smallskip
 Secondly, we set $M_h$ the maximum of the function $(u_1,u_2)\mapsto |||\GI_3-\GR^{'}_h(u_1,u_2)|||$ into the closed set $\big\{(u_1,u_2)\in \omega\; |\; h\le dist\big((u_1,u_2),\gamma_0\big)\le 2h\big\}$, and let $(s_1,s_2)$ be in this closed subset of $\omega$ such that
$$M_h= |||\GI_3-\GR^{'}_h(s_1,s_2)|||.$$ Applying the PoincarŽ-Wirtinger's inequality in the ball $B\big((s_1,s_2), 4h\big)$ we deduce that
$$\forall(s^{'}_1,s^{'}_2)\in B\big((s_1,s_2), h\big),\qquad |||\GR^{'}_h(s^{'}_1,s^{'}_2)-\GR^{'}_h(s_1,s_2)|||\le C||\nabla \GR^{'}||_{(L^2(B((s_1,s_2), 4h)))^3}.$$ The constant depends only on the PoincarŽ-Wirtinger's constant for a ball. 

\noindent  If $M_h$ is larger than $C||\nabla \GR^{'}||_{(L^2(B((s_1,s_2), 4h)))^3}$ we have
$$\eqalign{
&\pi h^2\big(M_h-C||\nabla \GR^{'}||_{(L^2(B((s_1,s_2), 4h)))^3}\big)^2\le 
||\GR^{'}_h-\GI_3|||^2_{(L^2(B((s_1,s_2), h)))^3}\cr
\le &||\GR^{'}_h-\GI_3||^2_{(L^2(\omega_{3h,\gamma_0}))^{3\times 3}}\le Ch^2||\nabla \GR^{'}||^2_{(L^2(\omega_{6h,\gamma_0}))^{3\times 3}}\cr}$$  then, in all the cases  we obtain 
$$M_h\le C||\nabla \GR^{'}||_{(L^2(\omega_{6h,\gamma_0}))^{3\times 3}}.$$ The constant does not depend on $h$ and $\GR^{'}$. The above inequalities show that there exists $h^{''}_0$  such that for any $h\le h^{''}_0$
$$|||\GR^{'}_h(s_1,s_2)-\GI_3|||\le C||\nabla \GR^{'}||_{(L^2(\omega_{6h,\gamma_0}))^{3\times 3}}\le {1/8}\qquad \hbox{for any } (s_1,s_2)\in \omega \quad \hbox{such that } \;\; h\le dist\big((s_1,s_2),\gamma_0\big)\le 2h.$$ By definition of $\GR_h$, that gives $|||\GR_h(s_1,s_2)-\GI_3|||\le {1/8}$.
\smallskip
Finally, for any $h\le \max(h^{'}_0,h^{''}_0)$ and for any  $(s_1,s_2)\in \overline{\omega}$ we have
$$dist(\GR_h(s_1,s_2) , SO(3))\le {1/8}.$$
Using (C.5) and (C.6) we obtain  (recall that $\big\|\cdot\big\|_2$ is the euclidian norm in $\R^3$)
$$\forall (s_1,s_2)\in\omega,\qquad \Big\|{\partial{\cal V}^{'}_h\over \partial s_\alpha}(s_1,s_2)-\GR^{'}_h(s_1,s_2)\Gt_\alpha(s_1,s_2)\Big\|_2 \le Ch||\phi||_{(W^{2,\infty}(\omega))^3}+C\big(||{\cal V}||_{(H^2(\omega_{3h}))^{3}}+||\GR^{'}||_{(H^1(\omega_{3h}))^{3\times 3}}\big)$$ where $\omega_{3h}=\{(s_1,s_2)\in \R^2\;|\; dist\big((s_1,s_2),\partial \omega\big)\le 3h\}$.

\noindent  We have
$${\partial{\cal V}_h\over \partial s_\alpha}-\GR_h\Gt_\alpha=(1-\psi_h)\Big({\partial{\cal V}^{'}_h\over \partial s_\alpha}-\GR^{'}_h \Gt_\alpha\Big)+{\partial\psi_h\over \partial s_\alpha}\big(\phi-{\cal V}^{'}_h\big).$$ Thanks to the above inequality, (C.4) and    again the estimate of  $|||\GR^{'}_h-\GI_3|||$ in the edge strip $h\le dist\big((s_1,s_2),\gamma_0\big)\le 2h$ we obtain for all $(s_1,s_2)\in\omega$
$$\eqalign{
&\Big\|{\partial{\cal V}_h\over \partial s_\alpha}(s_1,s_2)-\GR_h(s_1,s_2)\Gt_\alpha(s_1,s_2)\Big\|_2\cr
\le & C \big(h||\phi||_{(W^{2,\infty}(\omega))^3}+||{\cal V}||_{(H^2(\omega_{3h}))^{3}}+||\GR^{'}||_{(H^1(\omega_{3h}))^{3\times 3}}+||\phi-{\cal V}||_{(H^2(\omega_{5h,\gamma_0}))^{3\times 3}}\big).\cr}$$ 
The same argument as above    imply that there exists $h_0\le \max(h^{'}_0,h^{''}_0)$ such that for any $0<h\le h_0$ and for any $(s_1,s_2)\in \omega$ we have
$$\Big\|{\partial{\cal V}_h\over \partial s_\alpha}(s_1,s_2)-\GR_h(s_1,s_2)\Gt_\alpha(s_1,s_2)\Big\|_2 \le{1\over 8}.\leqno(C.11)$$  
From (C.4), (C.5), (C.6) and (C.7) there exists a positive constant $C$ which does not depend on $h$ such that
$$||\GR_h||_{(W^{1,\infty}(\omega))^{3\times 3}}\le {C\over h}\big\{||{\cal V}||_{(H^2(\omega))^3}+||\GR||_{(H^1(\omega))^{3\times 3}}\big\}.\leqno(C.12)$$ Now we can choose $h$ in term of $\delta$. We set
$$h=\theta\delta,\qquad \delta\in (0,\delta_0]$$ and we fixed $\theta$ in order to have $h\le h_0$ and to obtain the right hand side in (C.12) less than $\displaystyle {1\over 4\sqrt{2}c^{'}_1\delta}$ ( $c^{'}_1$ is given by (2.3)). It is well-known that there exists  a sequence $\big(\overline{V}_\delta\big)_{\delta\in (0,\delta_0]}$ such that $({\cal V}_\delta,\GR_\delta,\overline{V}_\delta)\in \D_{\delta,\gamma_0}$ and satisfying  the convergences in (C.1) and the estimate in (C.3).\fin

 \bigskip
\noindent{\Ggras References}
\medskip
\noindent [1] J.M. Ball, Convexity conditions and existence theorems in nonlinear elasticity. Arch. Ration. Mech. Anal. 63 (1976) 337-403.

\noindent [2] D. Blanchard, A. Gaudiello, G. Griso.  Junction of a periodic
family of elastic rods with a $3d$ plate. I. J. Math. Pures Appl. (9) 88 (2007), no 1, 149-190.

\noindent [3]   D. Blanchard, A. Gaudiello, G. Griso. Junction of a periodic
family of elastic rods with a thin plate. II. J. Math. Pures Appl. (9) 88 (2007), no 2, 1-33.

\noindent [4]  D. Blanchard,  G. Griso. Microscopic effects in the homogenization of the junction of rods and a thin plate.  Asympt. Anal. 56 (2008), no 1, 1-36.

\noindent [5] D. Blanchard, G. Griso. Decomposition of  deformations of thin rods. Application to nonlinear elasticity,  Ana. Appl. 7  (1) (2009) 21-71.

\noindent [6]  D. Blanchard, G. Griso.   Decomposition  of the deformations of a thin shell. Asymptotic behavior of the Green-St Venant's strain tensor. Journal of Elasticity: Volume 101 (2), (2010), 179-205.

\noindent [7]  D. Blanchard, G. Griso.   Justification of a simplified model for shells in nonlinear elasticity. C. R. Acad. Sci. Paris, Ser. I  348 (2010) 461-465.

\noindent [8] P.G. Ciarlet, Mathematical Elasticity, Vol. I, North-Holland, Amsterdam (1988).

\noindent [9] P.G. Ciarlet, Mathematical Elasticity, Vol. II. Theory of plates. North-Holland, Amsterdam (1997).

\noindent [10] P.G. Ciarlet, Mathematical Elasticity, Vol. III. Theory of shells. North-Holland, Amsterdam (2000).

\noindent [11] P.G. Ciarlet, Un modle bi-dimentionnel non linŽaire de coques analogue ˆ celui de W.T. Koiter, C. R. Acad. Sci. Paris, SŽr. I, 331  (2000), 405-410.

\noindent [12] P.G. Ciarlet, L. Gratie , C. Mardare. A nonlinear Korn inequality on a surface. J. Math. Pures Appl. 2006; 85: 2-16.

\noindent [13] P.G. Ciarlet and C. Mardare, An introduction to shell theory, Differential geometry: theory and applications, 94--184, Ser. Contemp. Appl. Math. CAM, 9, Higher Ed. Press, Beijing, 2008.

\noindent [14] G. Friesecke, R. D. James and S. MŸller. A theorem on geometric rigidity and the derivation of nonlinear plate theory from the three-dimensional elasticity. Communications on Pure and Applied Mathematics, Vol. LV, 1461-1506 (2002).

\noindent [15] G. Friesecke, R. D. James and S. MŸller, A hierarchy of plate models derived from nonlinear elasticity by $\Gamma$-convergence. (2005)

\noindent [16] G. Friesecke, R. D. James, M.G. Mora and S. MŸller, Derivation of nonlinear bending theory for shells from three-dimensionnal nonlinear elasticity by Gamma convergence, C. R. Acad. Sci. Paris, Ser. I 336 (2003).

\noindent [17] G. Griso. Asymptotic behavior of curved rods by the unfolding method. Math. Meth. Appl. Sci. 2004; 27: 2081-2110.

\noindent [18]  G. Griso. Asymptotic behavior of structures made of plates. Anal.  Appl. 3 (2005), 4, 325-356.

\noindent  [19]   G.  Griso. D\'ecomposition des d\'eplacements d'une poutre : simplification d'un probl\`eme d'\'elasticit\'e. C. R. Acad. Sci. Paris, Ser. II  333 (2005), 475-480.

\noindent [20]  G. Griso. Obtention d'Žquations de plaques par la mŽthode d'Žclatement appliquŽe aux Žquations tridimensionnelles. C. R. Acad. Sci. Paris, Ser. I 343 (2006), 361--366.

\noindent [21] G. Griso. Decomposition of displacements of thin structures.  J. Math. Pures Appl.  89 (2008) 199-233.

\noindent [22] H. Le Dret and A. Raoult, The nonlinear membrane model as variational limit of nonlinear three-dimensional elasticity. J. Math. Pures Appl. 75 (1995) 551-580.

\noindent [23] H. Le Dret and A. Raoult, The quasiconvex envelope of the Saint Venant-Kirchhoff  stored energy function. Proc. R. Soc. Edin., A 125 (1995) 1179-1192.

\noindent [24] J.E. Marsden and T.J.R. Hughes, Mathematical Foundations  of Elasticity, Prentice-Hall, Englewood Cliffs, (1983).

\noindent [25] O. Pantz,  On the justification of the nonlinear inextensional plate model. C. R. Acad. Sci. Paris SŽr. I Math. 332 (2001), no. 6, 587--592.

\noindent [26] O. Pantz,  On the justification of the nonlinear inextensional plate model. Arch. Ration. Mech. Anal. 167 (2003), no. 3, 179--209.

\bye